\documentclass[12pt]{article}
\usepackage{fullpage}
\usepackage{amsfonts}
\usepackage{amssymb}
\usepackage{hyperref}
\usepackage{bm}

\newtheorem{theorem}{Theorem}[section]
\newtheorem{lemma}[theorem]{Lemma}
\newtheorem{proposition}[theorem]{Proposition}
\newtheorem{remark}[theorem]{Remark}

\title{Arnold Diffusion in a priory chaotic Hamiltonian systems}
\author{V.Gelfreich$^*$ and D.Turaev$^{**}$\footnote{This work was supported by the grant RSF 14-41-00044.}\\
$^*$ Mathematics Institute, University of Warwick
\\
v.gelfreich@warwick.ac.uk \\[12pt]
$^{**}$
Department of Mathematics, Imperial College, London\\
d.turaev@imperial.ac.uk}
\begin{document}
\maketitle
\begin{abstract}
We assume that a symplectic real-analytic map has an invariant normally hyperbolic
cylinder and an associated transverse homoclinic cylinder. It is well known that
such cylinder is preserved under small perturbations. We prove that
for a generic real-analytic perturbation of this map the boundaries
of the cylinder are connected by trajectories of the perturbed map.
\end{abstract}


\newpage

\section{Introduction}

A Hamiltonian dynamical system is defined with the help of a Hamilton function $H:M\to\mathbb R$
on a symplectic manifold $M$ of dimension $2n$. Let $M_c$ be a connected component
of a level set  $\{H=c\}$. Since the Hamiltonian remains
constant along the trajectories of the Hamiltonian system, the set $M_c$ is invariant.
Depending on the Hamilton function $H$ and the energy $c$,
the restriction of the dynamics onto $M_c$ may vary from uniformly hyperbolic
(e.g., in the case of a geodesic flow on a surface of negative curvature) to completely
integrable.

Since Poincare's works it is accepted that a typical Hamiltonian systems does not have any additional integral of motion independent
of $H$ (unless the system posseses some symmetries and Noether theorem applies). On the other hand a generic Hamiltonian
system is nearly integrable in a neighbourhood of a totally elliptic equilibrium (a generic minimum or maximum of $H$) or totally
elliptic periodic orbit. Therefore, by the Kolmogorov-Arnold-Moser (KAM) theory, one should not expect from a given Hamiltonian system
to be ergodic (with respect to the Liouville measure) on every energy level \cite{MM1974}.
Indeed, the KAM-theory establishes that under sufficiently mild conditions a nearly integrable system possesses a set of invariant tori which has positive measure.

Recall that each of the KAM tori has dimension $n$ and, for $n>2$, does not divide $M_c$ which has dimension $(2n-1)$.
Moreover, for $n>2$ the complement to the union of all KAM tori is connected and dense in $M_c$.
So the KAM theory does not contradict to the existence of a dense orbit in $M_c$.
It is unknown whether such orbits really exist in nearly integrable systems.
The question goes back to Fermi \cite{Fermi1923} who suggested the following notion:
a Hamiltonian system is called {\em quasi-ergodic} if in every  $M_c$ any two open sets are connected by a trajectory.
This property is equivalent to topological transitivity of the Hamiltonian flow on $M_c$.
This property can also be restated in different terms: (a) in every $M_c$ there is a dense orbit or (b)
in every $M_c$ dense orbits form a residual subset.

Fermi conjectured \cite{Fermi1923} that quasi-ergodicity is a generic property of Hamiltonian systems,
but proved a weaker statement only: if a Hamiltonian system with $n>2$ degrees of freedom has the form
\begin{equation}\label{Eq:H}
H=H_0(I)+\varepsilon H_1(I,\varphi,\varepsilon),
\end{equation}
where $H_0$ is integrable and $(I,\varphi)$ are action-angle variables, then generically $M_c$ does not contain
an invariant $(2n-2)$-dimensional hyper-surface which is analytic in $\varepsilon$. The existence of such surface
would, obviously, prevent the quasi-ergodicity. However, it is not known whether a non-analytic invariant hyper-surface
can exist generically, so Fermi's quasi-ergodic hypothesis remains unproved. The recent papers
\cite{Mather2012,Kaloshin2012p,Cheng2013,Kaloshin2014,Marco2014} make an important
step in understanding of the underlying dynamics by showing that
in a generic (in a certain smooth category) case with $2\frac12$ or more degrees of freedom,
there are trajectories which visit an a-priory prescribed sequence of balls.
A paper \cite{KaloshinS2012} provides examples of an orbit
dense in a set of maximal Hausdorff dimension (equal to 5 as $n=3$).

This problem is closely related to the problem of stability of a totally elliptic fixed point
of a symplectic diffeomorphism, or stability of a periodic totally elliptic
periodic orbit for a Hamiltonian flow. It was proved in \cite{DL1983,D1988}
that stability can be broken by an arbitrarily small smooth perturbation.
It is believed that a totally elliptic periodic orbit is generically unstable
but the time scales for this instability to manifest itself are extremely long.

\medskip

The unperturbed system (\ref{Eq:H}) is described by the Hamiltonian $H=H_0(I)$.
Then the actions $I$ are constant along trajectories, so the equation $I=I_0$
defines an invariant torus, and the angles $\varphi$
are quasi-periodic functions of time with the frequency vector
$\omega_0(I)=H_0'(I)$.
KAM theory implies that the majority of invariant tori survive under perturbation.
Tori with rationally dependent frequencies are called {\em resonant} and are destroyed
by a typical perturbation. The frequency of a resonant torus satisfies a
condition of the form $\omega_0(I)\cdot\mathbf k=0$ for some $\mathbf k\in\mathbb Z^n\setminus\{0\}$.
The resonant tori form a ``resonant web", typically (e.g. if $\omega_0$ is a local
diffeomorphism) a dense set of measure zero.

Arnold's example \cite{Arnold} showed that a trajectory of the perturbed system  (\ref{Eq:H})
can drift along a resonance. This paper inspired a large number of studies in long-time stability of actions,
the problem which is known as ``Arnold diffusion''. It has been attracting a lot of attention recently
and we refer the reader to surveys in \cite{BK2005,DLS2006b,DLS2008b}
for a more detailed discussion.

It should be noted that the motion along the resonant web is very slow:
Nekhoroshev theory \cite{Nekhoroshev} provides a lower bound on the instability times.
Let $\{\cdot,\cdot\}$ denote the Poisson brackets. Then $\dot I=\{H,I\}=\varepsilon\{H_1,I\}$ is of the order of $\varepsilon$.
On the other hand, if the system satisfies assumptions of the KAM theory,
$|I(t)-I(0)|$ remains small for all times and the majority of initial conditions,
i.e., for the set of initial conditions of asymptotically full measure.
If $H$ satisfies assumptions of the Nekhoroshev theory, there are some exponents
$a,b>0$ such that $|I(t)-I(0)|<\varepsilon^{a}$ for all $|t|<\exp\varepsilon^{-b}$ and for
all initial conditions. This estimate establishes an exponentially large
lower bound for the times of Arnold Diffusion.

It is important to stress that the upper bound on the speed of Arnold diffusion strongly
depends  on the smoothness of the system. Indeed, the stability times are  exponentially large
in $\varepsilon^{-1}$ for analytic systems, but only polynomial bounds can be obtained in
the $C^k$ category. In particular, the papers \cite{Mather2012,Kaloshin2012p,Cheng2013}
study the Arnold diffusion for non-analytic Hamiltonians and therefore
the bounds established by the analytical Nekhoroshev theory are likely to be violated.
The problem of genericity of Arnold diffusion in analytic category remains
fully open. We believe the methods proposed in our paper will help to achieve
an advancement in the analytic case.

\medskip

In a neighbourhood of a simple resonance, normal form theory suggests existence
of a normally hyperbolic cylinder with a pendulum-like separatrix.
Bernard proved the existence of normally-hyperbolic cylinders in a priori stable
Hamiltonian systems the size of which is bounded from below independently
of the size of the perturbation \cite{Bernard2010a}.

A model for this situation is usually obtained by assuming that the integrable part
of the Hamiltonian already possesses a normally-hyperbolic cylinder and an associated homoclinic
loop (e.g. by considering $H_0=P(p,q)+h_0(I)$
where $P$ is a Hamiltonian of a pendulum). A system of this type is called
{\em a-priory unstable}. The drift of orbits along the cylinder has been actively
studied in the last decade
\cite{B2010,BounemouraP2012,CY2004,CY2009,DH2009a,DH2011,DLS2006b,T2004,Treschev2012},
including the problem of genericity of this phenomenon and instability times.
It should be noted that the Arnold diffusion is much faster in this case.

In these studies, a drifting trajectory typically stays most of the time near the
normally-hyperbolic cylinder, occasionally making a trip near a homoclinic loop.
The process can be described using the notion of a scattering map introduced
by Delshams, de la Llave and Seara in \cite{DLS2008a}. Earlier Moeckel \cite{M2002}
suggested that Arnold diffusion can be modeled by random application of
two area-preserving maps on a cylinder. In this way the deterministic
Hamiltonian dynamics is modeled by an iterated function system,
and the obstacles to a drift along the cylinder appear in the
form of common invariant curves \cite{M2002,L2007,NP2012}.

\bigskip

In our paper we depart from the near-integrable setting and study the dynamics
of an arbitrary exact symplectic map in a homoclinic channel, a neighbourhood of a normally-hyperbolic two-dimensional cylinder $A$
along with a sequence of homoclinic cylinders $B$ at the transverse intersection
of the stable and unstable manifolds of $A$. We conduct a rigorous reduction of the problem to
the study of an iterated function system and show that the existence of a drifting trajectory (i.e. the instability of the
Arnold diffusion type) is guaranteed when the exact symplectic maps of the cylinder $A$ that constitute
the iterated function system do not have a common invariant curve. The reduction scheme is in the same spirit as
in \cite{NP2012,LTG} (while the setting and proofs are different). The completely novel result is that the existence of
drifting orbits is a generic phenomenon, i.e. it holds for an open and dense subset of a neighbourhood, in the space of analytic symplectic
maps, of the given map with a homoclinic channel, provided the restriction of the map on the cylinder $A$ has a twist property.
All the known similar genericity results for the Arnold diffusion have been proven so far in the smooth category
and use the non-analiticity of the perturbations in an essential way.

In one respect, the situation we consider is more general than in the near-integrable setting, as we do not assume the
existence of a large set of KAM curves on the invaraint cylinder $A$. On the other side, as one can extract from the example
of \cite{DGT}, our assumption of the strong transversality of the homoclinic intersections which we need
in order to define the scattering maps that form the iteration function system seems to fail for a generic analytic
near-integrable system in a neighbourhood of a resonance in the a priori stable case. Therefore, our results do not
admit an immediate translation to the a priori stable case.

Rather, the problem we consider here is related to the {\em a-priory chaotic} case.
The term refers to a Hamiltonian of the form $H=H_0(p,q,I)+H_1(p,q,I,\theta)$, where
$p$ is the group of variables conjugate to the variables $q$, and $\theta$ is the angular variable conjugate to $I$.
Here $H_1$ is assumed to be small, and $H_0$ has the additional integral $I$ but is not completely integrable.
In particular, we speak about the a priori chaotic situation if
there exists a non-trivial hyperbolic set $\Lambda$ in the $(p,q,)$-space for each value
of the integral $I$ from some interval $[I_a,I_b]$. Then, adding the perturbation $H_1$ will typically cause the ``diffusion''
of $I$.

Namely, the existence of a non-trivial hyperbolic set $\Lambda$ is equivalent \cite{Shilnikov67} to the existence of a hyperbolic periodic
orbit whose stable and unstable manifolds have a transverse intersection along a homoclinic orbit. These two
orbits depend continuously on $I$, so we have a continuous family of periodic orbits and a corresponding continuous family of homoclinic orbits,
and we assume these two families are defined for all $I\in [I_a,I_b]$. A saddle periodic orbit in the $(p,q)$-space corresponds to a
saddle two-dimensional invariant torus in the full $(p,q,I,\theta)$-space. We take a cross-section of the form $F(p,q,I)=0$ to this torus; then
a Poincare map is defined on the cross-section, and the intersection of the torus with the cross-section is a saddle invariant circle of this map.
The circle is given by an equation $(p,q)=(p_I,q_I)=const$. The existence of a transverse homoclinic intersection of the stable and unstable
manifolds of the saddle periodic orbit in the $(p,q)$-space implies the existence of a transverse homoclinic intersection of the local stable
manifold of the invariant circle on the cross-section with a piece of its global unstable manifold. Since the Hamiltonian $H_0$ is
$\theta$-independent, i.e. the problem is symmetric with respect to rotation in $\theta$, each connected component of the homoclinic intersection is
a circle $(p,q)=(p_I^h,q_I^h)=const$. The union of the invariant circles $(p,q)=(p_I,q_I)$ over all $I\in[I_a,I_b]$
forms a normally-hyperbolic cylinder $A$ of the Poincare map; the union of the homoclinic circles $(p,q)=(p_I^h,q_I^h)=const$ over all
$I\in[I_a,I_b]$ forms a homoclinic cylinder $B$. It is easy to check that the {\em strong transversality conditions} from our Main Theorem
(Theorem~\ref{Thm:main}) are fulfilled by the homoclinic channel $\{A,B\}$.

Since $H_0$ is $\theta$-independent, and $I$ is the constant of motion, the Poincare map on the invariant cylinder $A$ has the form
$$\bar\theta =\theta+\omega(I), \qquad \bar I=I.$$
The map on the cylinder satisfies the twist condition
if $\omega'(I)\neq 0$ for all $I$ under consideration. Under this condition our Theorem~\ref{Thm:main}
implies the following result:\\
{\em Let $H_0$ be a real analytic function, e.g. it is holomorphic in a certain complex neighbourhood (analyticity domain) $Q$
of a bounded set in the real $(p,q,I,\theta)$ space that contains all the orbits with the initial conditions in the cylinders $A$ and $B$.
There exists $\varepsilon>0$ such that for an open and dense subset of the set of real analytic functions $H_1$ such that
$\sup_{cl(Q)} |H_1|<\varepsilon$ the system defined by the Hamiltonian $H_0+H_1$ has an orbit that goes from a small neighbourhood
of $I=I_a$ to a small neighbourhood of $I=I_b$.}

Since the unperturbed system preserves the value of $I$ and the normally-hyperbolic cylinder is an analytic manifold
foliated by the invariant circles, one can use the methods from the original Arnold paper \cite{Arnold} and construct
a special perturbation $H_1$ which does not change the system on the cylinder but creates heteroclinic connections
between different invariant curves on $A$. These connection can be used to form heteroclinic chains, and the existence
of drifting orbits follows. Our statement here is, in fact, much stronger, as it guarantees the existence of the drift
in $I$ for a {\em typical} analytic perturbation $H_1$.

A partial case is given by small periodic perturbations of autonomous Hamiltonian systems. Namely, consider a Hamiltonian of the form
$$
H=H_0(p,q)+H_1(p,q,t).
$$
We assume that the unperturbed Hamiltonian $H_0$ has, at every energy level from some interval $[I_a,I_b]$, a
hyperbolic periodic orbit with a transverse homoclinic; both orbits depend continuously on the value of the energy $I=H_0$.
The union of these periodic orbits over all $I$ is a normally-hyperbolic invariant cylinder $A$; the union of the homoclinics
provides a homoclinic cylinder $B$. The twist condition reads as $\partial T/\partial I \neq 0$ where $T(I)$ is the period of
the saddle periodic orbit in the energy level $H_0=I$. For all sufficiently small $H_1$ which are periodic in time $t$
the homoclinic channel $\{A,B\}$ persists and the twist condition holds. Then our theorem implies the existence of
a trajectory with the energy drifting from $I_a$ to $I_b$ for an open and dense (in the space of functions analytic in some
complex domain $Q$) set of small perturbations $H_1$. We note that the size of perturbation can be taken uniform for
all arbitrarily large frequencies of the perturbation, therefore the result holds true in the case where Melnikov type
computations are unavailable and the speed of the energy drift is bounded from above by the Nekhoroshev type estimates.

This problem is a natural extension of the Mather problem on the existence of trajectories
with unbounded energy in a periodically forced geodesic flow \cite{BT1999,DLS2000}. The criteria for the existence of
trajectories of the energy that grows up to infinity are known
\cite{BT1999,DLS2000,delaLlave2006,DLS2006a,GT2008,GTb} that are based on the assumption that the initial value of energy
is high enough. The above described result can be used for showing the generic existence of orbits of unbounded energy growth
for all possible initial energy values.

\smallskip

Other type of examples to our Main Theorem is given by small perturbations of direct products of two symplectic maps.
Thus, consider a 4-dimensional symplectic map which is a direct product of
a twist map and a standard map: $\Phi_0:(\varphi,I,x,y)\mapsto (\bar\varphi,\bar I,\bar x,\bar y)$:
\begin{equation}\label{Eq:4dsm}
\begin{array}{rclrcl}
\bar \varphi&=&\varphi+\omega(I),\qquad &\bar x&=&x+\bar y\,,\\
\bar I&=&I,\qquad&\bar y&=&y+k\sin x\,,
\end{array}
\end{equation}
where $k>0$ is a positive parameter and $\omega$ is an analytic function. We consider $\varphi$ and $x$
to be angular variables, so the map is a symplectic diffeomorphism of $(\mathbb T\times\mathbb R)^2$.
The map $\Phi_0$ preserves the value of the $I$ variable. Thus, the cylinder $A$ given by $x=y=0$ is invariant and
filled with invariant curves.

The situation is more challenging when the integrable twist map
is replaced by another standard map, so the new unperturbed map is given by $\Phi_0:(\varphi,I,x,y)\mapsto (\bar\varphi,\bar I,\bar x,\bar y)$:
\begin{equation}\label{Eq:4dsm-a}
\begin{array}{rclrcl}
\bar \varphi&=&\varphi+\bar I,\qquad &\bar x&=&x+\bar y\,,\\
\bar I&=&I+k_1\sin\varphi,\qquad&\bar y&=&y+k_2\sin x\,.
\end{array}
\end{equation}
The cylinder $A: x=y=0$ is still invariant but it is no longer filled with invariant curves. Instead
the cylinder contains a Cantor set of invariant curves provided $k_1$ is not too large.
These tori prevent trajectories of $\Phi_0$ from travelling in the direction of the $I$ axis.

\smallskip

The theory presented in this paper allows to treat both cases equally
and implies that an arbitrarily small generic analytic perturbation creates
trajectories which travel between regions $I<I_a$
and $I>I_b$ for any $I_a<I_b$ (provided $\omega'(I)$ is separated from $0$ for (\ref{Eq:4dsm}), and $k_2>C(4|k_1|+k_1^2)$ for (\ref{Eq:4dsm-a})).
Indeed, in order to apply Theorem~\ref{Thm:main} to
these examples, we notice, first, that the invariant cylinder $A$ is normally hyperbolic.
This cylinder has a stable and unstable separatrices $W^u(A)$ and $W^s(A)$
which coincide with the product of $A$ and the stable (reps., unstable) separatrix
of the standard map $W^{u,s}_{sm}$, so we can write (slightly abusing notation)
$W^s(A)=A\times W^s_{sm}$ and $W^u(A)=A\times W^u_{sm}$.
This product also describes the structure of the foliation
of $W^{u,s}(A)$ into strong stable and strong unstable manifolds of points in $A$.
For a point $v\in A$, we let $E^{uu}(v)=\{v\}\times W^u_{sm}$ and  $E^{ss}(v)=\{v\}\times W^s_{sm}$.
The assumption $k_2>C(4|k_1|+k_1^2)$ for (\ref{Eq:4dsm-a}) ensures that these
strong stable and strong unstable foliations remain $C^1$-smooth after the perturbation.

It is not too difficult to prove that the standard map has infinitely many transversal homoclinic
orbits for any $k>0$. Let $\bm p_h=(x_h,y_h)$ be one of these orbits.
The cylinder $B=A\times\{\, \bm p_h\,\}\subset W^u(A)\cap W^s(A)$ is homoclinic to $A$.
Since the strong stable and strong unstable foliations of a point $ v\in A$ coincide with the product of
the base point and the separatrices of the standard map, we see that $(v,p_h)\in E^{ss}(v)\cap E^{uu}(v)$,
and the cylinder $B$ satisfies the strong tansversality assumption described in the next Section
giving rise to a {\em simple} homoclinic intersection (defined in the Next section).
Then Theorem~\ref{Thm:main} implies that generic perturbation of $\Phi_0$ has orbits travelling in
the direction of the cylinder $A$.

Similar maps were considered in Easton et al. \cite{EMR2001} (motivated by the ``stochastic pump model" of Tennyson et al. \cite{TLL1979}).
In \cite{EMR2001} the existence of drift orbits was shown for all non-integrable Lagrangian perturbations provided $k_2$ is large enough
(i.e. in the ``anti-integrable'' limit). Our methods allow to obtain the drifting orbits without the large $k_2$ assumption, i.e.
without a detailed knowledge of the dynamics of the system. The price for that is that we do not have a constructive description of the
class of perturbation which lead to the Arnold diffusion instability, although we know that these perturbations form a dense (and open)
subset of the set of all analytic perturbations that keep the homoclinic channel intact.

\section{Set up, assumptions, and results}

Consider a real-analytic diffeomorphism $\Phi:\Sigma\to \mathbb R^{2d}$, $d\ge2$, defined on
an open set $\Sigma\subseteq\mathbb R^{2d}$. We assume that $\Phi$ preserves the standard symplectic form
$\Omega$, and that $\Phi$ is exact (the latter is always true if $\Sigma$ is e.g. simply-connected).
Let $\Phi$ have an invariant smooth two-dimensional
cylinder $A$ diffeomorphic to $\mathbb S^1\times [0,1]$. We denote the corresponding
embedding $\mathbb S^1\times [0,1]\to \Sigma$ as $\psi$. The cylinder $A$ is bounded
by a pair of invariant curves $\psi(\mathbb S^1\times {0})$ and $\psi(\mathbb S^1\times {1})$. We will use the notation
$\partial A$ for this pair of curves, and $int(A)$ for $A\backslash \partial A$. We also
denote $F_0=\Phi|_{A}$.

We assume that the cylinder $A$ is normally-hyperbolic. Namely, we assume that at each point
$v\in A$ the tangent space is decomposed into direct sum of three non-zero subspaces:
$\mathbb R^{2d}=N^c_v\oplus N^u_v\oplus N^s_v$, where $N^c_v$ is the two-dimensional
tangent to $A$ at the point $v$, while $N^{s,u}$ are certain subspaces which depend continuously on $v$,
and these subspaces fields are invariant with respect to the derivative $\Phi'$ of the map,
i.e. $\Phi'N^s_v=N^s_{F_0(v)}$, $\Phi'N^u_v=N^u_{F_0(v)}$ (the condition $\Phi'N^c_v=N^c_{F_0(v)}$ is
fulfilled automatically because of the invariance of $A$ with respect to $\Phi$).
We assume that for some choice of norms in $N^{s,u,c}$
there exist $\alpha>1$ and $\lambda\in(0,1)$ such that at every point $v\in A$
\begin{equation}\label{lf0}
\|F'_0(v)\|<\alpha, \qquad \|(F_0'(v))^{-1}\|<\alpha,
\end{equation}
\begin{equation}\label{ldl}
\|\Phi'|_{N^s_v}\|<\lambda, \qquad \|(\Phi'|_{N^u_v})^{-1}\|<\lambda,
\end{equation}
and
\begin{equation}\label{nhm}
\alpha^2\lambda<1.
\end{equation}
Note that these assumptions are stronger compared to the standard definition of a normally-hyperbolic
manifold. In particular, the large spectral gap condition (\ref{nhm}) implies
the $C^1$-regularity of the strong stable and strong unstable foliations
while the general case implies H\"older continuity only (see e.g. \cite{PSW2011}). The $C^1$-regularity means, in our case,
that the leaves of the foliation are smooth and, importantly, the field of tangents to the leaves is also smooth, which implies that
for any two cross-sections transverse to the foliation the map between the cross-sections which is defined by the leaves of the foliation
is a diffeomorphism.

We also note that, according to the above definition, $A$ is {\em symmetrically} normally-hyperbolic,
i.e., in (\ref{lf0}),(\ref{ldl}) the same pair of exponents $\alpha$, $\lambda$ bounds both $\Phi'$ and $(\Phi')^{-1}$.
This symmetric form of the spectral gap assumption is chosen because the map $\Phi$ is symplectic.
We will prove (Proposition \ref{smplca}) that the symmetric spectral gap implies that
the restriction of the symplectic form on $A$ is non-degenerate, i.e. $A$ is a symplectic submanifold of $\mathbb R^{2d}$,
therefore the map $F_0=\Phi|_{A}$ inherits the (exact) symplecticity of $\Phi$.

We are quite confident that the factors $\lambda$ and $\alpha$ may be allowed to depend on the point $v$ on the manifold.
However, for simplicity, we conduct the proofs only for the case of $v$-independent $\lambda$ and $\alpha$.

The points in a small neighbourhood of the normally-hyperbolic cylinder $A$,
whose forward iterations tend to $A$ without leaving the neighbourhood,
form a smooth (at least $C^2$ in our case) invariant manifold,
the local stable manifold $W^s_{loc}\supset A$, which is tangent to $N^s\oplus N^c$ at the points of $A$ (see e.g. \cite{HPS}).
The points whose backward iterations tend to $A$, without leaving the neighbourhood, form the local unstable
smooth invariant manifold $W^u_{loc}\supset A$, which is tangent to $N^u\oplus N^c$ at the points of $A$.
The invariant cylinder $A$ is the intersection of $W^u_{loc}$ and $W^s_{loc}$. In each of the manifolds
there exists a uniquely defined $C^1$-smooth invariant foliation transverse to $A$, the strong-stable invariant foliation
$E^{ss}$ in $W^s_{loc}$ and the strong-unstable invariant foliation $E^{uu}$ in $W^u_{loc}$, such that
for every point $v\in A$ there is a unique leaf of $E^{ss}_v$ and a unique leaf of $E^{uu}_v$ which pass through
this point and are tangent to $N^s_v$ and, respectively, $N^u_v$ (see \cite{SSTC}).

The global stable and unstable manifolds of $A$ are defined by iterating the
local invariant manifolds: $W^u(A):=\bigcup_{m\geq 0} \Phi^m W^u_{loc}$ and
$W^s(A):=\bigcup_{m\geq 0} \Phi^{-m} W^s_{loc}$.
As the global stable and unstable manifolds are obtained from the local ones by the iterations of $\Phi$,
the locally defined invariant foliations $E^{ss}$ and $E^{uu}$ are extended in a unique way on the whole
of $W^s(A)$ and, respectively, $W^u(A)$.

\bigskip

Let us now assume that {\em the symmetrically normally-hyperbolic cylinder $A$ has a homoclinic},
i.e.,  the intersection of $W^u(A)$ and $W^s(A)$ has a point $x$ outside $A$. If $W^u(A)$ and $W^s(A)$
are transverse at $x$, the implicit function theorem implies that $x$ has an open neighbourhood
$U_x$ in $W^u(A)\cap W^s(A)$, which is diffeomorphic to a two-dimensional disk.

Note that for any $x\in W^u(A)\cap W^s(A)$
there exists a unique leaf of $E^{uu}_x$ and a unique leaf of $E^{ss}_x$
which pass though this point. We will call the homoclinic intersection at $x$
{\em strongly transverse} if
\begin{equation}\label{sttc1}
\mathbb R^{2d}={\mathcal T}_xE^{ss}_x\oplus {\mathcal T}_xE^{uu}_x\oplus {\mathcal T}_x (W^u(A)\cap W^s(A)),
\end{equation}
i.e. the leaf $E^{uu}_x$ is transverse to $W^u(A)$ and the leaf $E^{ss}_x$ is transverse to $W^s(A)$ at
the point $x$. In other words, the foliation $E^{uu}$ is transverse to the small disc $U_x$ in $W^u(A)$ and
the foliation $E^{ss}$ is transverse to $U_x$ in $W^s(A)$. Thus, the strong transversality implies that
the holonomy maps $\pi^s: U_x\to A$ and $\pi^u: U_x\to A$ (projections by the leaves of the foliations $E^{ss}$ and
$E^{uu}$, respectively) are diffeomorphisms.

\medskip

Let $\bar A\subset int(A)$ be a compact invariant sub-cylinder in $A$, i.e. it is a closed
region in $int(A)$ bounded by two non-intersecting invariant essential curves $\gamma^+$ and $\gamma^-$.
We say that a smooth manifold $B\subset W^u(int(A))\cap W^s(int(A))\setminus A$ is a {\em homoclinic
cylinder, simple relative to the cylinder $\bar A$}, if the following assumptions hold:
\begin{itemize}
\item[{[S1]}]
The {\em strong transversality} condition (\ref{sttc1}) holds for all $x\in B$.
\item[{[S2]}]
For every point $x\in\bar A$, the corresponding leaf of the foliation $E^{uu}$ intersects the homoclinic cylinder $B$
at exactly one point each, and no two points in $B$ belong to the same leaf of the foliation $E^{ss}$.
In other words, the projection $\pi^{s}_B: B\to int(A)$ by the leaves of the strong-stable foliation is injective,
and the projection $(\pi^{u}_B)^{-1}: \bar A \to B$ by the leaves of the strong-unstable foliation is well-defined.
\item[{[S3]}] The image of $\bar A$ by the map $\pi^{s}_B \circ (\pi^{u}_B)^{-1}$ contains an
essential curve.
\end{itemize}

Condition [S1] implies that the projections $\pi^{s,u}_B: B\to int(A)$ are local diffeomorphisms. Condition
[S2] implies that the maps $(\pi^s_B)^{-1}$ and $(\pi^u_B)^{-1}$ are well-defined (single-valued).
Thus, following \cite{DLS2008a}, we define
the {\em scattering map}
$F_B=\pi^s_B \circ (\pi^u_B)^{-1}$, which is a diffeomorphism $\bar A\to F_B(\bar A)\subset int(A)$.
Condition [S3] implies that the image of any essential curve by the scattering map is also an essential curve, i.e.
the scattering map is homotopic to identity on $\bar A$. In fact, it is easy to see that conditions [S2] and [S3] are equivalent to the
requirement that the scattering map is a homotopic to identity diffeomorphism $\bar A\to F_B(\bar A)$.

In particular, it follows that the set $F_B(\bar A)$ is a cylinder bounded by two non-intersecting essential curves,
$F_B(\gamma^{\pm})$. We will also show (Proposition \ref{Prop:sympl}) that
the scattering map $F_B$ is an exact symplectomorphism. This implies that the cylinder $F_B(\bar A)$
has the same area as $\bar A$ and that the intersection of $F_B(\bar A)$ and $\bar A$ is non-empty.

We will call the orbit of $B$ by the map $\Phi$ a {\em homoclinic sequence  of cylinders}. Since the foliations $E^{ss}$
and $E^{uu}$ are invariant, it follows that each set $\Phi^n(B)$ in this sequence is also a simple homoclinic cylinder.
Indeed, the invariance of the foliations means that $\Phi(\pi^s(x))=\pi^s(\Phi(x))$ and
$\Phi(\pi^u(x))=\pi^u(\Phi(x))$ for every point $x$ from $W^s(A)$ or, respectively, $W^u(A)$. Hence,
\begin{equation}\label{psubb}
\pi^{s,u}_{\Phi(B)}=\Phi\circ \pi^{s,u}_B\circ \Phi^{-1}.
\end{equation}
Thus, the scattering map satisfies
\begin{equation}\label{sccon}
F_{\Phi(B)}=\Phi\circ F_B\circ\Phi^{-1}.
\end{equation}
As we see, the scattering maps corresponding to any two different cylinders in the same homoclinic sequence are smoothly conjugate to each other,
so the fulfilment of the simplicity conditions for one of the cylinders implies the fulfilment of the simplicity conditions for the other one.

\medskip

Note that the fulfilment of condition [S2]
depends both on the choice of the invariant cylinder $\bar A$ and on the choice of the ambient invariant cylinder $A$:
the latter must be large enough to incorporate $F_B(\bar A)$. We will need a stronger version of this requirement.
Namely, we assume that $F_B(\bar A)\subseteq \hat A$ where $\hat A\supseteq \bar A$ is a compact invariant sub-cylinder
of $A$. One can express this property as follows:
\begin{equation}\label{wsuab}
W^u(\bar A)\cap B \subset W^s(\hat A)\cap B.
\end{equation}
We will also require that the homoclinic cylinder $B$ is simple relative to the cylinder $\hat A$.

We further assume that one can choose
symplectic coordinates $(y,\varphi)$ on $A$ such that $F_0=\Phi|_A$ will have a twist property.
Birkhoff theorem \cite{Herman1983} implies then that the boundary curves $\gamma^{\pm}$
of the invariant sub-cylinder $\bar A$ are graphs of Lipschitz functions $y^{\pm}(\varphi)$.

\bigskip

We will study behaviour of a generic map from a small neighbourhood $\mathcal V$ of the map $\Phi$
in the space of real-analytic
exact symplectic maps. We define the topology in this space as follows. Take any compact neighbourhood $K$
of the cylinder $A$ in $\Sigma$, which includes the sequence of homoclinic cylinders,
and let $Q$ be any compact complex neighbourhood of $K$ such that
$\Phi$ is holomorphically extended onto some open neighbourhood of $Q$. The neighbourhood $\mathcal V$
of the map $\Phi$ in the space of real-analytic exact symplectomorphisms consists of all holomorphic maps
$cl(Q)\to\mathbb C^{2d}$ that are sufficiently close to $\Phi$ everywhere on $Q$ and are
real on $K$, preserve the standard symplectic form in $\mathbb R^{2d}$, and are exact; two maps
belonging to this neighbourhood are close if they are uniformly close on $Q$.

We will assume that the invariant cylinder $A$ persists for every map in $\mathcal V$.
Namely, we assume that every map from $\mathcal V$ has, in a small neighbourhood of $A$,
a pair of non-intersecting essential invariant
curves which depend continuously on the map (as continuous curves), and these curves form
the boundary $\partial A$ of the normally-hyperbolic cylinder $A$ for our original map $\Phi$.
One shows that for every map from $\mathcal V$ these curves bound a uniquely defined
invariant cylinder which depends continuously on the map. We will continue to denote this cylinder as $A$.
It is natural just to assume that $A$ is bounded by two KAM-curves;
then the required persistence of the boundary curves is guaranteed.

Note that the invariant sub-cylinder $\bar A$ does not
need to be preserved when the map $\Phi$ is perturbed, as the invariant curves
$\gamma^{\pm}$ that form its boundary can, in principle, be destroyed by a small perturbation
(we do not assume that they are KAM-curves). However, this is not important for us (neither we require the
preservation of the cylinder $\hat A$).
We just define the curves $\gamma^{\pm}\subset A$ in an arbitrary way
for all maps from $\mathcal V$ so that they will depend continuously on the map,
and for the map $\Phi$ they will bound the invariant cylinder $\bar A$.

We say that a map connects two open sets $U^-$ and $U^+$ if the image of $U^-$
by some forward iteration of the map intersects $U^+$.

\begin{theorem}[main theorem]\label{Thm:main}
Let $A$ be a symmetrically normally-hyperbolic compact invariant cylinder for the exact symplectic map $\Phi$,
and let $F_0=\Phi|_A$ have a twist property. Let $\bar A$ be a compact invariant sub-cylinder
$\bar A \subset int(A)$ bounded by two non-intersecting essential curves. Let
$B\subset W^u(int(A))\cap W^s(int(A))$ be a
homoclinic cylinder, simple relative to $\bar A$. Suppose $\bar A\subseteq \hat A\subseteq A$
where $\hat A$ is a compact invariant sub-cylinder such that $W^u(\bar A)\cap B \subset W^s(\hat A)\cap B$ and
$B$ is simple relative to $\hat A$. Suppose $A$ persists for every map $\tilde\Phi$ from
a sufficiently small neighbourhood $\mathcal V$ of $\Phi$ in the space of real-analytic
exact symplectomorphisms $\Sigma\to \mathbb R^{2d}$. Let $\gamma^-$ and $\gamma^+$
be simple essential curves in $A$, which depend continuously on the map $\tilde\Phi$
and which coincide with the boundary of $\bar A$ when $\tilde\Phi=\Phi$.
Then, for an open and dense subset $\tilde{\mathcal V}$ of $\mathcal V$, each map $\tilde\Phi\in \tilde{\mathcal V}$
connects every two open neighbourhoods $U^-$ of $\gamma^-$ and $U^+$ of $\gamma^+$.
\end{theorem}

\begin{remark} {\em It is obvious
that given any two open sets $U^+$ and $U^-$ the set of maps that connect $U^-$ and $U^+$
is open. We, however, show that intersection of all these sets (over all possible choices of the neighbourhoods
$U^-$ and $U^+$ of the given curves $\gamma^-$ and $\gamma^+$) is open and dense
in $\mathcal V$, i.e. given any map from $\mathcal V$ there exists an open set of arbitrarily small
perturbations of this map within $\mathcal V$ such that each of these perturbations creates,
for each pair of neighbourhoods $U^-$ and $U^+$ of the curves $\gamma^\pm$, an orbit that connects $U^-$ and $U^+$.}
\end{remark}

\begin{remark} {\em Statements similar to Theorem \ref{Thm:main} are known
for non-analytic (smooth) case, see e.g. \cite{CY2004,CY2009,NP2012}. The main difference between the analytic
and smooth case is that the class of perturbations small in the real-analytic sense is narrower than
the class of perturbations that are small in the $C^{\infty}$-sense (e.g. our theorem implies the similar statement in
the smooth category). Crucially, for a typical real-analytic map
the normally-hyperbolic invariant cylinder $A$ is not analytic (it has only finite smoothness), so
no real-analytic perturbations can vanish on $A$. This makes the methods that have been
used in the non-analytic case \cite{CY2004,CY2009,NP2012} inapplicable.}
\end{remark}

\begin{remark} {\em The symplectic diffeomorphism $\Phi$ can be a Poincare map for a Hamiltonian
flow inside a level of constant energy. The methods of this paper can be generalised to
show that  if the Poincare map $\Phi$ for some Hamiltonian system satisfies the assumptions of the theorem,
then a generic small perturbation of the Hamiltonian function $H$
in the space of real-analytic Hamiltonians leads to creation of connecting orbits.}
\end{remark}

\bigskip

The strategy of the proof of our main theorem is as follows. We show in Proposition~\ref{Prop:manyb}
that the existence of one homoclinic cylinder $B$ which is simple relative to the two invariant cylinders $\bar A$ and $\hat A$
such that (\ref{wsuab}) holds
implies the existence of infinitely many secondary homoclinic cylinders (all belonging to different homoclinic sequences)
which are simple relative to $\bar A$. Thus, we will no longer use the existence of the invariant cylinder $\hat A$ and
will further consider $N\geq 8$ homoclinic cylinders $B_1,\dots, B_N$, all of which are simple relative to the same
compact invariant cylinder $\bar A$ and all belong to different homoclinic sequences, i.e.
$\Phi^m (B_i)\cap B_j=\emptyset$ for all $m$ and all $i,j=1,\dots,N$ such that $i\neq j$.
Given each of the homoclinic cylinders $B_n$, simple relative to the compact invariant
sub-cylinder $\bar A\subset int(A)$, we define the corresponding scattering map
$F_n:\bar A\to A$. By condition [S1], it is a local diffeomorphism. By condition [S2]
this map is a bijection, hence it is a diffeomorphism of $\bar A$ onto the set $F_n(\bar A)$, which
is a cylinder bounded by a pair of essential curves as follows from condition [S3].
Obviously, condition [S1] implies that the scattering maps are, in fact, defined in some open neighbourhood
$A'$ of $\bar A$ in $A$.

\medskip

It is a standard fact from the theory of normal hyperbolicity \cite{HPS} that any strictly-invariant normally-hyperbolic
compact smooth manifold with a boundary can be extended to a locally-invariant normally-hyperbolic manifold without a boundary.
In our case this means that the smooth embedding $\psi$ that defines the ivariant cylinder $A=\psi(\mathbb S^1\times [0,1])$
can be extended onto $\mathbb S^1\times I$ where $I$ is an open interval containing $[0,1]$, and the image
$\tilde A=\psi(\mathbb S^1\times I)\supset A$ is normally-hyperbolic and locally-invariant with respect to the map $\Phi$.
Here, by the local invariance we mean that there exists a neighbourhood $Z$ of $\tilde A$ such that the iterations of each point
of $\tilde A$ stay in $\tilde A$ until they leave $Z$.  An important property of the locally-invariant normally-hyperbolic manifold
without a boundary is that it
persists at $C^2$-small perturbations, i.e. all maps $C^2$-close to $\Phi$ have a locally-invariant normally-hyperbolic
cylinder $\tilde A\subset Z$ which depends on the map continuously as a $C^2$-manifold\footnote{Throughout this paper we assume
the large spectral gap assumption (\ref{nhm})
in the notion of normal hyperbolicity. This guarantees the $C^2$-smoothness of the manifold, and the $C^1$-smoothness
of the corresponding strong-stable and strong-unstable invariant foliations for every map $C^2$-close to $\Phi$.}.
The continuous dependence on the map implies
that the cylinder $\tilde A$ remains symplectic and symmetrically normaly-hyperbolic for all maps $C^2$-close to $\Phi$, e.g. for all maps from
$\mathcal V$.

Note that the normal hyperbolicity implies that $\tilde A$ contains all the orbits that
never leave $Z$. In particular, any invariant curve that lies in $Z$ must lie in $\tilde A$. Thus, the pair
of simple non-intersecting essential invariant curves that, by our assumption, exist in $Z$ for all maps from $\mathcal V$ must lie in $\tilde A$.
Hence, for all maps from $\mathcal V$ these curves bound a compact invariant sub-cylinder $A\subset \tilde A$,
which is uniquely defined by the choice  of the invariant boundary curves and depends continuously on the map.
The stable
and unstable manifolds and the strong-stable and strong-unstable foliations of $A$ also depend continuously, in the $C^1$-topology, on the map.
The transversality condition [S1] implies that the $C^1$-smooth homoclinic cylinders $B_1,\dots, B_N$ also persist and depend continuously
on the map. Thus, the scattering maps $F_1,\dots, F_N$ are defined
on the same open subset $A'$ of $A$ for all maps from $\mathcal V$, and remain exact symplectomorphisms. Note that
$A'$ contains the region bounded by the curves $\gamma^{\pm}$.

\medskip

Take any map $\tilde\Phi\in \mathcal V$. Let
$(v_s)_{s=0}^{m}\subset A$ be an orbit of iterated function system
$\{F_0,\ldots,F_N\}$, i.e. for each $s=0,\dots, m-1$ there exists $n_s=0,\dots,N$ such that $v_{s+1}=F_{n_s}(v_s)$; if $n_s\neq 0$, then we
always require $v_s\in A'$, so the corresponding map $F_{n_s}$ is well-defined.
In Section~\ref{Se:shadowing} we show that for any such orbit and any $\varepsilon>0$, there is a
point $x_0$ and a positive integer $\ell$ such that
$$
\mathrm{dist}(x_0,v_0)<\varepsilon, \mbox{~~~~and~~~} \mathrm{dist}(\tilde\Phi^{\ell}(x_0),v_m)<\varepsilon
$$
(see Lemma \ref{lmsw}). Note that we do not use
hyperbolicity or index arguments in this lemma. We also do not use the symplecticity of the
maps $F_1,\dots, F_N$, nor the twist property of the map $F_0$. However, the fact that $F_0$ is
an area-preserving map of a bounded invariant domain $A$
is crucial here, as we use the Poincare Recurrence Theorem in an essential way (we first prove a certain week shadowing result,
Lemma \ref{shadprop}, that holds without any assumptions on the map $F_0$, then Lemma \ref{lmsw} is deduced from it
in the case of area-preserving $F_0$).

According to this shadowing lemma (Lemma \ref{lmsw}), in order to show that two open sets are connected by the map $\tilde\Phi$, it
is sufficient to show that the intersections of these sets with $A$ are connected by orbits of the iterated function system $\{\,F_0,\ldots,F_N\,\}$.
A generalisation (Theorem~\ref{Thm:Birkhoff}) of a classical Birkhoff theorem
states that if $F_n$, $n=0,\dots,N$ are homotopic to identity, exact symplectomorphisms, and $F_0$ is a twist map, then
for any two essential curves $\gamma^\pm\subset A'$ there is
a trajectory of the iterated function system with $v_0\in\gamma^-$ and $v_m\in\gamma^+$
unless the functions $F_n$ have a common invariant essential curve.

Thus, if the maps $F_0,\dots, F_N$ have no common invariant essential curves between $\gamma^-$ and $\gamma^+$, every pair
of neighbourhoods, $U^-$ of $\gamma^-$ and $U^+$ of $\gamma^+$, is connected by the map $\tilde\Phi$. Theorem \ref{Thm:Birkhoff}
also implies that this property of the absence of a common invariant essential curve is open.

Theorem~\ref{Thm:nocurve} establishes that this property is also dense in $\mathcal V$ (provided $N\geq 8$).
Thus, for every map $\tilde\Phi$ from an open and dense subset of $\mathcal V$, the corresponding scattering maps $F_1,\dots, F_N$ ($N\geq 8$)
and $F_0$ do not have any common essential invariant curve. As we just explained,
this implies that every two neighbourhoods $U^\pm$ of $\gamma^\pm$ are
connected by each such map $\tilde\Phi$, and Theorem~\ref{Thm:main} follows.

\medskip

Theorem~\ref{Thm:nocurve} is the crucial step in the proof of main theorem~\ref{Thm:main}.
An analogue of Theorem~\ref{Thm:nocurve} for generic non-analytic maps can be derived from \cite{CY2004,CY2009,NP2012}.
However, the methods of those papers cannot be used in the analytic case (as the real-analytic perturbations cannot,
in general, vanish on the finitely smooth normally-hyperbolic cylinder). Therefore, we
develop a completely different perturbation technique in order to prove Theorem~\ref{Thm:nocurve} for the analytic case.

\section{Estimates in a neighbourhood of a symmetrically normally-hyperbolic invariant cylinder}

\subsection{Fenichel coordinates, cross form of the map, and estimates for the local dynamics}

We start with the analysis of the behaviour in a small neighbourhood of a normally-hyperbolic cylinder.
We do not need analyticity or symplecticity of the map in this and the next Sections.

Let $A$ be a compact, symmetrically normally-hyperbolic, smooth, invariant cylinder of a $C^r$-smooth map $\Phi$ ($r\geq 2$).
As we mentioned, $A$ can be extended to a larger, smooth
normally-hyperbolic locally-invariant cylinder $\tilde A$. Let us introduce coordinates in a small neighbourhood of $A$
such that this larger invariant cylinder is straightened; moreover, the local stable and unstable manifolds $W^{s,u}_{loc}$
are straightened as well, along with the strong-stable and strong-unstable foliations $E^{ss}$ and $E^{uu}$ on them.
Note that the foliations are at least $C^1$. Thus, the straightening of the
foliations means that one can introduce $C^1$-coordinates $(u,v,z)$
in a neighbourhood of $A$ such that the manifold $W^s_{loc}$
will have equation $z=0$, the manifold $W^u_{loc}$ will be given by $u=0$, and the leaves of the foliations
$E^{ss}$ and $E^{uu}$ will all have the form $\{z=0,v=const\}$ and, respectively,
$\{u=0,v=const\}$ (cf. \cite{JT2009}). The cylinder $A$ thus lies in $(u=0,z=0)$.
Here $v=(\varphi,y)$ with $\varphi\in \mathbb S^1$ being the angular variable,
and $y$ taking values from some (maybe $\varphi$-dependent) interval $I$ of the real line.

Note that the manifolds $W^u(A)$ and $W^s(A)$ can be non-orientable, so when we say we introduce global coordinates $(u,v,z)$
in a neighbourhood of $A$, we mean that $\varphi=0$ and $\varphi=2\pi$ are glued by means of some linear involution in the space
of coordinates $(v,z)$. This does not affect any estimates below.

\medskip

Since the manifolds and foliations under consideration are invariant with respect to the map $\Phi$,
it follows that in our coordinates the map near $A$ has the form
\begin{equation}\label{tmp}
\bar u=h_1(u,z,v),\qquad \bar z=h_2(u,z,v),\qquad
\bar v=F_0(v)+h_3(u,z,v),
\end{equation}
where $h_{1,2,3}$ and $F_0$ are $C^1$-functions such that
\begin{equation}\label{gg}
\begin{array}{c}
h_1(0,z,v)\equiv 0,\qquad h_2(u,0,v)\equiv 0,\\
h_3(0,z,v)\equiv 0,\qquad h_3(u,0,v)\equiv 0.
\end{array}
\end{equation}
It is seen from (\ref{tmp}),(\ref{gg}) that the manifolds $u=0$ and $z=0$, and the foliations
$(v=const, u=0)$ and $(v=const, z=0)$ are indeed invariant with respect to the map. Note that
by construction
$$\left.\frac{\partial h_1}{\partial u}\right|_{z=0}=\Phi'|_{N^s},\qquad
\left.\frac{\partial h_2}{\partial z}\right|_{u=0}=\Phi'|_{N^u},$$
therefore, we obtain from (\ref{ldl}) that (for an appropriate choice of norms)
$$
\left\|\frac{\partial h_1}{\partial u}\right\|<\lambda,
\qquad \left\|\left(\frac{\partial h_2}{\partial z}\right)^{-1}\right\|<\lambda.
$$

Using the implicit function theorem for small $u$ and $z$,
one may resolve the $\bar z$-equation in (\ref{tmp}) with respect to $z$.
Therefore, the map $\Phi:(u,v,z)\mapsto(\bar u,\bar v,\bar z)$ can be written, in some neighbourhood of
the closed invariant cylinder $A$ in the following cross form:
\begin{equation}\label{mpc}
\bar u=p(u,v,\bar z),\qquad z=q(u,v,\bar z),
\end{equation}
\begin{equation}\label{ggf}
\bar v=F_0(v)+f(u,v,\bar z),
\end{equation}
such that
\begin{equation}\label{ff}
p(0,v,\bar z)\equiv 0,\qquad q(u,v,0)\equiv 0,
\end{equation}
\begin{equation}\label{ggf2}
f(0,v,\bar z)\equiv 0,\qquad f(u,v,0)\equiv 0,
\end{equation}
\begin{equation}\label{lf00}
\|F'_0(v)\|<\alpha, \qquad \|(F_0'(v))^{-1}\|<\alpha,
\end{equation}
\begin{equation}\label{lgl}
\left\|\frac{\partial p}{\partial u}\right\|<\lambda,\qquad \left\|\frac{\partial q}{\partial \bar z}\right\|<\lambda,
\end{equation}
where
\begin{equation}\label{alla}
\alpha\lambda<1, \qquad 0<\lambda<1<\alpha,
\end{equation}
see (\ref{lf0})-(\ref{nhm}). Let $Z_\delta$ denote a $\delta$-neighbourhood of $A$.

\begin{lemma}\label{Le:1}
There is $\delta_0>0$ such that for any $\delta\in(0,\delta_0)$ and any $k\geq 0$
the following statements hold.
\begin{enumerate}
\item
Any trajectory of length $k$ such that $(u_i,v_i,z_i):=\Phi^i(u_0,v_0,z_0)\in Z_\delta$ for $i=0,\dots,k$
satisfies the following estimates for $i=0,\dots,k$:
\begin{equation}\label{lamb}
\|u_i\|\leq \delta \lambda^i, \qquad \|z_i\|\leq \delta \lambda^{k-i},
\end{equation}
\begin{equation}\label{shad0}
\|v_i-F_0^iv_0\|\leq\delta (\alpha\lambda)^{k/2}, \qquad \|v_i-F_0^{i-k}v_k\|\leq\delta (\alpha\lambda)^{k/2}.
\end{equation}
\item
The orbit $(u_i,v_i,z_i)$ is determined in a unique way for any given $u_0$, $v_0$, $z_k$ such that
$\|u_0,z_k\|\leq \delta$ and $v_0\in A$, as well as for any
given $u_0$, $v_k$, $z_k$ such that $\|x_0,z_k\|\leq \delta$ and $v_k\in A$.
\item
Moreover, as $k\to+\infty$,
\begin{equation}\label{dery}
\left\|\frac{\partial z_0}{\partial (u_0,v_0)}\right\|+
\left\|\frac{\partial (u_k,v_k)}{\partial z_k}\right\|
\to 0,\quad
\left\|\frac{\partial u_k}{\partial (v_k,z_k)}\right\|+
\left\|\frac{\partial (v_0,z_0)}{\partial u_0}\right\|\to 0,
\end{equation}
uniformly for all $\|u_0,z_k\|\leq\delta$ and all $v_k\in A$ or $v_0\in A$.
\item
We also have for all $k$ large enough
\begin{equation}\label{derz}
\left\|\frac{\partial z_0}{\partial z_k}\right\|\leq \lambda^k,\qquad
\left\|\frac{\partial (u_k,v_k)}{\partial (u_0,v_0)}\right\|\leq \alpha^k
\end{equation}
(at any given $(u_0,v_0)$ in the first inequality, and at any given $z_k$ in the second one), and
\begin{equation}\label{derv}
\left\|\frac{\partial u_k}{\partial u_0}\right\|\leq \lambda^k,\qquad
\left\|\frac{\partial (v_0,z_0)}{\partial (v_k,z_k)}\right\|\leq \alpha^k,
\end{equation}
(at any given $(v_k,z_k)$ in the first inequality, and at any given $u_0$ in the second one).
\end{enumerate}
\end{lemma}

\noindent {\em Proof.} Using (\ref{mpc}), we get
\begin{equation}\label{Eq:traj}
u_{i+1}=p(u_i,v_i,z_{i+1}),\qquad z_i=q(u_i,v_i,z_{i+1}),
\qquad
v_{i+1}=F_0(v_i)+f(u_i,v_i,z_{i+1}),
\end{equation}
for all $i=0,\dots,k-1$. Equations  (\ref{ff}) and (\ref{lgl}) imply that
\begin{equation}\label{pqd}
\|u_{i+1}\|=\|p(u_i,v_i,z_{i+1})\|\leq \lambda \|u_i\|, \qquad \|z_{i}\|=\|q(u_i,v_i,z_{i+1})\|\leq \lambda \|z_{i+1}\|\,.
\end{equation}
Since $\|u_0\|,\|z_k\|\leq \delta$, it follows that the orbit $\{(u_i,z_i,v_i)\}_{i=0}^k$ satisfies (\ref{lamb}).

\medskip

For the future convenience let us define
\begin{equation}\label{Eq:C0}\!\!\!\!\!\!\!\!
C_0(\delta)=\max\left\{\sup_{Z_\delta}\|p'_v\|,\ \sup_{Z_\delta}\|p'_{\bar z}\|,\
\sup_{Z_\delta}\|q'_u\|,\ \sup_{Z_\delta}\|q'_{\bar z}\|,\
\sup_{Z_\delta}\|f'_u\|,\ \sup_{Z_\delta}\|f'_v\|,\
\sup_{Z_\delta}\|f'_{\bar z}\|\right\},
\end{equation}
and note that $C_0(\delta)$ can be made as small as we need by decreasing $\delta$
because (\ref{ff}) and (\ref{ggf2}) imply that
$p'_v=0$, $p'_{\bar z}=0$, $q'_u=0$, $q'_v=0$, $f'_u=0$, $f'_v=0$, $f'_{\bar z}=0$ at $(u=0, z=0)$, for all $v\in A$.

In order to establish inequalities (\ref{shad0}), let $V_i:=v_i-F_0^iv_0$. In particular $V_0=0$.
Equation (\ref{Eq:traj}) implies
\begin{equation}\label{uuii}
\|V_{i+1}\|\leq \sup_{v\in A} \|F_0'(v)\| \cdot \|V_i\|+\|f(u_i,v_i,z_{i+1})\|.
\end{equation}
Then equation (\ref{ggf2})  implies
\begin{eqnarray*}
\|f(u_i,v_i,z_{i+1})\|&\leq& \sup_{(u,v,z)\in Z_\delta}  \|f'_u\| \cdot\|u_i\|\le C_0(\delta)\|u_i\|\,,
\\
\|f(u_i,v_i,z_{i+1})\|&\leq &\sup_{(u,v,z)\in Z_\delta} \|f'_z\|\cdot\|z_{i+1}\|\le C_0(\delta)\|z_{i+1}\|\,.
\end{eqnarray*}
Thus, by (\ref{lamb}),
\begin{equation}\label{fmin}
\|f(u_i,v_i,z_{i+1})\|\leq \delta\ C_0(\delta)\ \min\{\lambda^i,\lambda^{k-i-1}\}\,.
\end{equation}
Now, by (\ref{lf00}) and (\ref{fmin}), we may rewrite (\ref{uuii}) as
$$\|V_{i+1}\|\leq \alpha \|V_i\|+\delta\ C_0(\delta)\ \min\{\lambda^i,\lambda^{k-i-1}\}.$$
Since $V_0=0$, we find (see (\ref{alla})) that for all $1\leq j\leq k$
\begin{eqnarray*}\!\!\!\!
\|V_j\|&\leq& \delta\ C_0(\delta) \sum_{0\leq i\leq j-1} \alpha^{j-i-1}\min\{\lambda^i,\lambda^{k-i-1}\}\leq
\delta C_0(\delta)\sum_{0\leq i\leq k-1} \alpha^{k-i-1}\min\{\lambda^i,\lambda^{k-i-1}\}=\\
&&=\delta\ C_0(\delta) \left\{\sum_{~~0\leq i\leq (k-1)/2} \!\!(\alpha\lambda)^{k-i-1}+
\alpha^{k-1}\sum_{(k-1)/2< i\leq k-1}\!\!(\lambda/\alpha)^i\right\}\ \leq \ \delta(\alpha\lambda)^{k/2}\,,
\end{eqnarray*}
when $\delta_0$ is chosen small enough to ensure $\displaystyle
\frac{C_0(\delta)}{\sqrt{\alpha\lambda}}\left[\frac{1}{1-\alpha\lambda}+\frac{\lambda}{\alpha-\lambda}\right]\leq 1$.
The first of inequalities (\ref{shad0}) is proved.
The second inequality follows immediately by
the symmetry of the problem (if we replace the map $\Phi$ by its inverse, then $F_0$ changes to $F_0^{-1}$, $i$ to $(k-i)$,
$(u_0,z_k)$ to $(z_k,u_0)$ and $v_0$ to $v_k$).

\bigskip

Given $u_0$, $v_0$, $z_k$, the orbit $\{(u_i,z_i,v_i)\}_{i=0}^k$ is
a fixed point of the operator $$Q: \{(u_i,v_i,z_i)\}_{i=0}^k\mapsto \{(\hat u_i,\hat v_i,\hat z_i)\}_{i=0}^k,$$
which acts on a sequence $\{(u_i,v_i,z_i)\}_{i=0}^k$ by
\begin{equation}\label{tmpa}
\left\{\begin{array}{l}
\hat u_{i+1}=p(u_i,v_i,z_{i+1}),\qquad
\hat z_i=q(u_i,v_i,z_{i+1}),\\ [4pt]
\hat v_{i+1}=F_0(v_i)+ f(u_i,v_i,z_{i+1})
\qquad\qquad\qquad\qquad\mbox{for $i=0,\dots,k-1$},\\[4pt]
\hat u_0=u_0,\qquad \hat v_0=v_0,\qquad \hat z_k=z_k.
\end{array}\right.
\end{equation}

Recall that $v=(y,\varphi)$, where $\varphi\in \mathbb S^1$, and $y$ runs an interval $I$ such that
for all sufficiently small $\delta$ the points in the $\delta$-neighbourhood $Z_\delta$ of the
cylinder $A$ have the $y$-coordinates strictly inside $I$.
It will, however, be convenient to extend the functions $p,q,F_0,f$ in (\ref{mpc}),(\ref{ggf})
to all $y\in \mathbb R^1$ in such a way that they will remain smooth, will have uniformly continuous derivatives,
identities (\ref{ff}),(\ref{ggf2}) will hold, and estimates (\ref{lf00}),(\ref{lgl}) will stay true with a margin of
safety. We assume that this extension is done without changing the functions $p,q,F_0,f$ inside the interval $I$
of $y$-values. Thus, if for a fixed point
of the operator $Q$ given by (\ref{tmpa}) with the functions $p,q,F_0,f$, which are defined now for all $y$,
every $(u_i,v_i,z_i)$, $i=0,\dots,k$, belongs to $Z_\delta$, then this
fixed point is an orbit of the original map $\Phi$.

It is also convenient to consider the lift of the original map so that
$\varphi$ runs the whole real axis and the functions $p$, $q$ and $F_0+f-v$ are periodic
in $\varphi$. So, in the analysis of the operator $Q$ given by (\ref{tmpa}), we assume $v\in \mathbb R^2$.

\medskip

Denote by $X=X_{k,u_0,v_0,z_k}$ the set of all sequences $\{(u_i,v_i,z_i)\}_{i=0}^k$ with the given
value of $(u_0,v_0,z_k)$, which also satisfy $\|u_i,z_i\|\leq \delta$ for all $i=0,\dots,k$. By (\ref{pqd}),
if $\|u_i,z_i\|\leq \delta$ for all $i=0,\dots,k$, then $\|\hat u_i,\hat z_i\|\leq \delta$ for all $i=0,\dots,k$
as well, thus $QX\subseteq X$.
Let us show that the operator $Q$ is contracting on $X$ in the norm
$$\|\{(u_i,v_i,z_i)\}_{i=0}^k\|_\alpha=\max_{i=0,\dots,k} \alpha^{-i} \|u_i,v_i,z_i\|.$$
Indeed, in this norm
\begin{eqnarray*}
\|Q'\|_\alpha&\leq&
\max\Biggl\{
\alpha^{-1}\left\|\frac{\partial p}{\partial u}\right\|
+
\alpha^{-1}\left\|\frac{\partial p}{\partial v}\right\|
+
\left\|\frac{\partial p}{\partial \bar z}\right\|,\qquad
\left\|\frac{\partial q}{\partial u}\right\|
+
\left\|\frac{\partial q}{\partial v}\right\|
+
\alpha\left\|\frac{\partial q}{\partial \bar  z}\right\|,\\ &&\qquad\qquad
\alpha^{-1}\left\| F_0'\right\|+
\alpha^{-1}\left\|\frac{\partial f}{\partial u}\right\|
+
\alpha^{-1}\left\|\frac{\partial f}{\partial v}\right\|
+
\left\|\frac{\partial f}{\partial  \bar z}\right\|
\Biggr\}\leq
\\
&\le&
\max\Bigl\{
\alpha^{-1}\lambda
+
\alpha^{-1}C_0(\delta)
+C_0(\delta),\
2C_0(\delta)
+
\alpha\lambda
,\
\alpha^{-1}\left\| F_0'\right\|+
\alpha^{-1}
C_0(\delta)
+
2C_0(\delta)
\Bigr\},
\end{eqnarray*}
where, for the derivatives in the right-hand side, we use the supremum norm taken over all $(u,v,\bar z)$ such that $\|u,\bar z\|\leq \delta$,
and $C_0(\delta)$ is defined by (\ref{Eq:C0}).
By (\ref{ff})-(\ref{alla}), if $\delta$ is sufficiently small,
then $\|Q'\|_\alpha<1$ uniformly for every element from $X$, independently of the value of $k\geq 0$.
Since the set $X$  is convex, it follows that the operator $Q$ is indeed contracting.

Thus, by contraction mapping principle, given any $(u_0,v_0,z_k)$ such that $\|u_0,z_k\|\leq\delta$
there exists indeed a unique length-$k$ orbit with the given values of $u_0$, $v_0$ and $z_k$.
We already proved that this orbit must satisfy (\ref{lamb}) and (\ref{shad0}). Since $v_0\in A$ implies
$F_0^iv_0\in A$ for all $i=0,\dots,k$ by the invariance of $A$ with respect to $F_0$, estimates (\ref{lamb}) and (\ref{shad0})
imply that the orbit lies in $Z_\delta$ as required.

By the symmetry of the problem, given any $(u_0,v_k,z_k)$ such that $\|u_0,z_k\|\leq\delta$ and $v_k\in A$,
there, as well, exists a unique length-$k$ orbit with the given values of $u_0$, $v_k$ and $z_k$,
and this orbit lies in $Z_\delta$.

As a fixed point of a smooth contracting operator, the obtained orbit must depend smoothly on all data
on which the operator depends smoothly. So, $(u_i,v_i,z_i)$ depend smoothly on $(u_0,v_0,z_k)$
(and, by the symmetry of the problem, on $(u_0,v_k,z_k)$ as well). To complete the proof
of the lemma, it remains to prove estimates
(\ref{dery}),(\ref{derz}) and (\ref{derv}).

\medskip

We will prove only the first limit in (\ref{dery}), as the second one follows then due to the symmetry of the problem
with respect to change of $\Phi$ to $\Phi^{-1}$; it is also enough to prove only (\ref{derz}), as (\ref{derv})
follows then by the symmetry.
Denote $\beta_i=\|\partial (u_i,v_i)/\partial (u_0,v_0)\|$, $\gamma_i=\|\partial z_i/\partial (u_0,v_0)\|$,
where the derivatives are taken at $z_k$ fixed. By differentiating (\ref{Eq:traj}),
we obtain
$$
\beta_{i+1}\leq \left\|\frac{\partial (p,F_0+f)}{\partial(u,v)}\right\|\ \beta_i\ +
\left\|\frac{\partial (p,f)}{\partial \bar z}\right\|\ \gamma_{i+1},\qquad
\gamma_{i}\leq \left\|\frac{\partial q}{\partial \bar z}\right\|\ \gamma_{i+1}+
\left\|\frac{\partial q}{\partial (u,v)}\right\|\ \beta_i,
$$
where the derivatives are taken at $(u,v,\bar z)=(u_i,v_i,z_{i+1})$. Since $u_i$ and $z_i$ satisfy (\ref{lamb}),
we obtain from (\ref{ff})-(\ref{lgl}) that for sufficiently small $\delta$
(independent of $i$ and $k$)
\begin{equation}\label{bgg}
\beta_{i+1}\leq (\alpha-\rho) \beta_i + \mu_i \gamma_{i+1},\qquad
\gamma_{i}\leq (\lambda-\rho) \gamma_{i+1}+ \mu_{k-i-1} \beta_i,
\end{equation}
where $\rho$ is a small positive constant, and
\begin{equation}\label{phips}
\mu_j=\sup_{\|u\|\leq \delta\lambda^j,\; (u,v,z)\in Z_\delta}
\left\|\frac{\partial (p,f)}{\partial \bar z}\right\|+
\sup_{\|z\|\leq \delta\lambda^j,\; (u,v,z)\in Z_\delta}
\left\|\frac{\partial q}{\partial (u,v)}\right\|\;;
\end{equation}
it follows from (\ref{ff}),(\ref{ggf2}) that
\begin{equation}\label{dimp}
\mu_j\to0 \mbox{~~~as~~} j\to +\infty.
\end{equation}
Recall also that, by definition,
\begin{equation}\label{bncd}
\beta_0=1,\qquad \gamma_k=0.
\end{equation}

Define the sequence $M_j$ by the rule
\begin{equation}\label{mnrl}
M_{j+1}=\alpha\lambda M_j +\mu_j,
\end{equation}
for an arbitrarily chosen $M_0$. As $\alpha\lambda<1$, it follows from (\ref{dimp}) that
\begin{equation}\label{dimn}
M_j\to0 \mbox{~~~as~~} j\to +\infty.
\end{equation}
By (\ref{bgg})
$$\gamma_{i}-M_{k-i}\beta_i\leq \frac{\lambda-\rho}{1- \mu_i M_{k-i-1}} (\gamma_{i+1}-M_{k-i-1}\beta_{i+1})+
\left[\mu_{k-i-1}-M_{k-i}+\alpha\frac{\lambda-\rho}{1- \mu_i M_{k-i-1}} M_{k-i-1}\right]\beta_i.
$$
As both the sequences $\mu_j$ and $M_j$ tend to zero, it follows that
\begin{equation}\label{mmlim}
\lim_{k\to+\infty} \max_{i=0,\dots,k-1} \mu_i M_{k-i-1}=0,
\end{equation}
so if $k$ is large enough, then $\mu_i M_{k-i-1}<\rho/\lambda<1$ for all $i=0,\dots,k-1$. Thus,
$$\gamma_{i}-M_{k-i}\beta_i\leq \lambda (\gamma_{i+1}-M_{k-i-1}\beta_{i+1})+
\left[\mu_{k-i-1}-M_{k-i}+\alpha\lambda M_{k-i-1}\right]\beta_i,
$$
which, by (\ref{mnrl}), implies
$$\gamma_{i}-M_{k-i}\beta_i\leq \lambda (\gamma_{i+1}-M_{k-i-1}\beta_{i+1}),$$
hence, for all $k$ large enough, for every $i=0,\dots,k-1$
\begin{equation}\label{kgbi}
\gamma_i-M_{k-i}\beta_i\leq \lambda^{k-i} (\gamma_k-M_0\beta_k),
\end{equation}
in particular
\begin{equation}\label{kgb}
\gamma_0-M_k\beta_0\leq \lambda^k (\gamma_k-M_0\beta_k).
\end{equation}
Now, by (\ref{bncd}), we have $\gamma_0\leq M_k$, so (\ref{dimn}) implies
$\partial z_0/\partial (u_0,v_0)\to 0$ as $k\to+\infty$, which agrees with (\ref{dery}). Note also that
by (\ref{kgbi}) we have $\gamma_{i+1}\leq M_{k-i-1}\beta_{i+1}$. By (\ref{bgg}),(\ref{mmlim}), this
gives us that for all $k$ large enough, for every $i=0,\dots,k-1$
$$\beta_{i+1}\leq\alpha\beta_i,$$
which (see (\ref{bncd})) implies the second inequality in (\ref{derz}).

It remains to estimate $\partial(u_k,v_k,z_0)/\partial z_k$ as $k\to+\infty$. To this aim,
we denote now $\beta_i=\|\partial (u_i,v_i)/\partial z_k\|$,
$\gamma_i=\|\partial z_i/\partial z_k\|$,
where the derivatives are taken at $(u_0,v_0)$ fixed. Then by differentiating (\ref{mpc}),(\ref{ggf}), we
will obtain exactly the same inequalities (\ref{phips}) as before, hence the estimate (\ref{kgb})
holds at all sufficiently large $k$ for the newly defined $\gamma_i,\beta_i$.
However, instead of (\ref{bncd}) we have now
$$\beta_0=0,\qquad \gamma_k=1.$$
Thus, we find from (\ref{kgb}) that
$$\gamma_0\leq\lambda^k, \qquad \beta_k\leq 1/M_0$$
for all $k$ sufficiently large. This immediately gives us the first inequality in (\ref{derz}), and since $M_0$ can be
taken arbitrary, we also obtain that $\partial(u_k,v_k)/\partial z_k\to 0$ as $k\to+\infty$, which finishes the proof of (\ref{dery}).
\hfill $\Box$

\subsection{``Lambda-lemma''}
The following analogue of the ``lambda-lemma'' \cite{Palis,Wig} is extracted from Lemma~\ref{Le:1}.

\begin{proposition}\label{Lemma:lambda}
Given any surface $L$ of the form $u=w(v,z)$, where $w$ is a smooth function
defined for all $v\in A$ and all small $z$, the images of $L$ by the map $\Phi$ converge to $W^u_{loc}(A)$
in the $C^1$-topology. If a surface $L$ has the form $z=w(v,u)$, where $w$ is a smooth function
defined for all $v\in A$ and all small $u$, then the images $\Phi^{-m}L$ converge (in $C^1$) to $W^s_{loc}(A)$
as $m\to+\infty$
\end{proposition}

\noindent
{\em Proof.}  By the symmetry of the problem, it is enough to consider only the case where $L$ is a surface of the form
$u=w(v,z)$. By Lemma~\ref{Le:1}, given any $(u_0,v_k,z_k)$ the corresponding orbit
$(u_i,v_i,z_i)$ is defined uniquely. Denote as $\eta_k$
the operator that sends $(u_0,v_k,z_k)$ to $(v_0,z_0)$, and as $\xi_k$ the operator that sends
$(u_0,v_k,z_k)$ to $u_k$. The point $(u,v,z)$ belongs to $\Phi^k L$ if and only if
$u_0=w(v_0,z_0)$, i.e. the equation of $\Phi^k L$ is
\begin{equation}\label{lkph}
u_k=\xi_k(u_0,v_k,z_k)
\end{equation}
where $u_0$ is defined from
\begin{equation}\label{lamblem}
u_0=w(\eta_k(u_0,v_k,z_k)).
\end{equation}
By (\ref{lamb}) and (\ref{dery}),
$$\|\eta_k\|+\|\partial \eta_k/\partial u_0\|\to 0 \mbox{~~as~~} k\to+\infty,$$
therefore at each $k$ large enough equation (\ref{lamblem}) defines $u_0$ uniquely as a smooth
function of $(v_k,z_k)$. It follows from (\ref{derv}) that
$$\|\partial u_0/\partial (v_k,z_k)\|=O(\alpha^k).$$
Thus, equation (\ref{lkph}) defines $u_k$ as a smooth function $w_k(v_k,z_k)$, for all $\|z_k\|\leq\delta$
and $v_k\in A$. By (\ref{lamb}), $\|u_k\|\to0$ as $k\to+\infty$. Moreover, since by (\ref{derv}) and (\ref{dery})
we have $\|\partial \xi_k/\partial u_0\|=O(\lambda^k)$ and $\|\partial \xi_k/\partial (v_k,z_k)\|\to 0$ as $k\to+\infty$,
it follows that
$$\left\|\frac{dw_k}{d (v_k,z_k)}\right\|\leq \left\|\frac{\partial \xi_k}{\partial u_0}\right\|\cdot
\left\|\frac{\partial u_0}{\partial (v_k,z_k)}\right\|+
\left\|\frac{\partial \xi_k}{\partial (v_k,z_k)}\right\|=
O((\alpha\lambda)^k)+ \left\|\frac{\partial \xi_k}{\partial (v_k,z_k)}\right\|\to 0$$
as $k\to+\infty$ (recall that $\alpha\lambda<1$). We see that for all $k$ large enough the surface $\Phi^k L$ is
given by the equation $u=w_k(v,z)$ where $w_k$ tends to zero along with the first derivative as $k\to +\infty$.
Since equation of $W^u_{loc}$ is $u=0$, this proves the proposition. ~~$\sqcap\!\!\!\!\sqcup$

\subsection{Secondary homoclinic cylinders}

Proposition \ref{Lemma:lambda} allows to establish a sufficient condition for the existence of infinitely many homoclinic cylinders.
Let $\bar A\subset A$ and $\hat A\subset A$ be a pair of compact invariant cylinders, $\bar A\subset \hat A$. Let the intersection
of $W^u(A)$ and $W^s(A)$ be non-empty
and contain a homoclinic cylinder $B$ which is simple relative to $\hat A$.
Since the manifolds $W^u(A)$ and $W^s(A)$ are invariant with respect to $\Phi$,
all iterations $\Phi^m B$ also belong to the intersection $W^u(A)\cap W^s(A)$, and this homoclinic sequence of cylinders converges to $A$
as $m\rightarrow \pm\infty$. As we explained before (see (\ref{sccon})), if $B$ is a simple homoclinic cylinder relative to $\bar A$
then the cylinders $\Phi^m B$, for any $m$, also has this property.

\begin{proposition}\label{Prop:manyb}
Let condition (\ref{wsuab}) hold for the pair of cylinders $\bar A$ and $\hat A$.
Then there are infinitely many homoclinic cylinders $B_i$,
each corresponds to a simple (relative to $\bar A$) intersection of $W^u(A)$ with $W^s(A)$, and none of the cylinders belongs
to the homoclinic sequence corresponding to another cylinder: $B_i\cap \Phi^m B_j=\emptyset$ for every $m$ and every $i\neq j$.
\end{proposition}

\noindent{\em Proof.}
Take the homoclinic cylinder $B$ and consider its iterations $B^-=\Phi^{-m_-}(B)$ and $B^+=\Phi^{m_+} (B)$ ($m_{\pm}>0$) which lie
in a small neighbourhood of $A$. Thus, $B^-\subset W^u_{loc}$ and $B^+\subset W^s_{loc}$.
By the relative to $\hat A$ simplicity of $B$, the cylinder $B^-$ intersects transversely the leaves of $E^{uu}$ in $W^u_{loc}(\hat A)$ at
a single point each. Since $W^u_{loc}$ has the equation $u=0$ and the leaves of
the foliation $E^{uu}$ in $W^s_{loc}$ are given by $\{u=0,v=const\}$,
it follows that there is a piece $W$ of the manifold $W^s(A)$ which contains the homoclinic cylinder
$B^-$ and has the form $z=w(v,u)$ where $w$ is a smooth
function defined for all $v$ from some neighbourhood of $\hat A$ and all small $u$.
By Proposition~\ref{Lemma:lambda} (where the invariant cylinder $A$ is replaced by the invariant cylinder $\hat A$),
the images $W_i=\Phi^{-i} W$ by the backward iterations of $\Phi$ accumulate on $W^s_{loc}(\hat A)$ in $C^1$.
By (\ref{wsuab}), $W^s_{loc}(\hat A)\cap B^+\supset W^u(\bar A)\cap B^+$. Therefore,
each of $W_i$ with $i$ sufficiently large intersects $W^u(\bar A)$ near $B^+$ transversely. Since $W_i$
are, by construction, pieces of $W^u(A)$, this gives
us the sought infinite set of homoclinic cylinders $B_i$ (converging to the cylinder $B^+$). Since $W_i$ are $C^1$-close
to $W^s_{loc}(\hat A)$ near $B^+$, it follows from the relative to $\hat A$ symplicity of $B^+$ that
$W_i$ intersect transversely each leaf of the foliation $E^{uu}$ in $W^u(\bar A)$, the uniqueness of the
intersections is also inherited.

Thus, the scattering maps $F_i: \bar A\to int(A)$ are defined for each of the cylinders $B_i$. In order to check
the simplicity of the homoclinic intersection at $B_i$, we need to show that the projections $\pi^{s}_{B_i}: B_i\to int(A)$ by
the leaves of the strong-stable foliation are injective for all $i$ (condition [S2]), and that the scattering maps
are homotopic to identity (condition [S3]). To check the injectivity, notice that
\begin{equation}\label{psiubbb}
\pi^{s}_{\Phi^i(B_i)}=\Phi^i\circ \pi^{s}_{B_i}\circ \Phi^{-i}
\end{equation}
by (\ref{psubb}). So, it is enough to show the injectivity of $\pi^{s}_{\Phi^i(B_i)}$. To do this, note that
the cylinders $\Phi^i(B_i)$ are close to $B^-$ at large $i$, so the maps $\pi^{s}_{\Phi^i(B_i)}$ are close to $\pi^s_{B^-}$,
and the latter map is injective by the simplicity of the homoclinic intersection at $B$.

It remains to show that the scattering maps $F_i$ are homotopic to identity. As we just mentioned, the maps
$\hat\pi_i^s=\pi^s_{\Phi^i(B_i)}\circ (\pi^s_{B^-})^{-1}$ are close to identity at large $i$. The same is true for the maps
$\hat\pi_i^u=\pi^u_{B^+}\circ(\pi^u_{B_i})^{-1}$. Using (\ref{psiubbb}), we find
\begin{equation}\label{lngfr}
F_i=\pi^s_{B_i}\circ (\pi^u_{B_i})^{-1}=\Phi^i\circ \hat\pi_i^s  \circ F_{B^-}\circ \pi^u_{B^-}\circ \Phi^{-i}\circ
(\pi^s_{B^+})^{-1} \circ F_{B^+}\circ \hat\pi^u_i,
\end{equation}
where $F_{B^+}=\pi^s_{B^+}\circ (\pi^u_{B^+})^{-1}$ and
$F_{B^-}=\pi^s_{B^-}\circ (\pi^u_{B^-})^{-1}$ are the scattering maps corresponding to the cylinders $B^+$ and $B^-$.
By the simplicity of the homoclinic intersection at $B$, these maps are homotopic to identity diffeomorphisms. The map
$\Phi^i$ in formula (\ref{lngfr}) acts in a small neigbourhood $Z$ of $A$ and is homotopic to identity in $Z$. The maps
$\pi^s_{B^+}$ and $\pi^u_{B^-}$ are projections along the foliations in the local stable and unstable manifolds, so they
are homotopic to identity in $Z$. Thus, all the maps in the right-hand side of formula (\ref{lngfr}) are homotopic to identity,
which implies that the scattering maps $F_i$ are homotopic to identity for all $i$ large enough. This proves the claim.
\hfill$\Box$

This proposition shows that the assumptions of Theorem \ref{Thm:main} imply the existence of an infinite series of homoclinic
cylinders which are simple relative to $\bar A$. For our purposes, the existence of $N\geq 8$ such cylinders is enough, so it will be
our standing assumption for the rest of the paper. We do not need the auxilliary invariant cylinder $\hat A$ anymore.

\section{Shadowing in the homoclinic channel\label{Se:shadowing}}

\subsection{Homoclinic channel}

Let  $B_1,\dots,B_N$ be homoclinic cylinders, each corresponds to a simple homoclinic intersection relative to
the compact invariant subcylinder $\bar A$, and none of the cylinders $B_n$ belongs to the homoclinic sequence corresponding
to another cylinder. Let us repeat the definition of the scattering maps $F_n$.
Since the homoclinic intersections are simple, it follows that
two diffeomorphisms, $\pi^{u}_n$ and $\pi^{s}_n$, from $B_n$ into $A$
are defined for every $n$ by the leaves of the foliations $E^{uu}$ and $E^{ss}$, respectively. Namely, $v=\pi^{u}_n(x)$ if the points
$x\in B_n$ and $v\in A$ belong to the same leaf of the foliation $E^{uu}$, and
$v=\pi^{s}_n(x)$ if $x\in B_n$ and $v\in A$ belong to the same leaf of the foliation $E^{ss}$
(the smoothness of the maps $\pi^{s}_n$ and $\pi^{u}_n$ and their inverse maps follows from the transversality of
the intersections of the leaves with $B_n$). By assumption,
$\bar A\subset \pi^{s}_n(B_n)$. Thus,
for each homoclinic cylinder $B_n$ we have a diffeomorphism $F_n=\pi^{s}_n\circ (\pi^{u}_n)^{-1}$
which acts from $\bar A$ into $A$. In fact, as the strong transversality condition [S1] is open,
there is a neighbourhood $A'$ of $\bar A$ such that all the scattering maps $F_1,\dots, F_N$ are diffeomorphisms of $A'$ into $A$.

Take sufficiently large positive $m_+$ and $m_-$ such that all the cylinders $B^+_n=\Phi^{m_+}(B_n)$ and $B^-_n=\Phi^{-m_-}(B_n)$
($n=1,\dots,N$) lie in the $\delta$-neighbourhood of $A$, where $\delta$ is small enough.
As $B^+_n\in W^s_{loc}$ and $B^-_n\in W^u_{loc}$, it follows that $z=0$ on $B^+_n$, and $u=0$ on $B^-_n$.
Since the homoclinics are simple, the cylinder $B^+_n$ intersects the leaves $\{v=const\}$ of the foliation $E^{ss}$
in $W^s_{loc}$ transversely, no more than at one point each, hence $B^+_n$ is given by
$B^+_n: \{u=u^+_n(v), z=0\}$,
where $u^+$ is a smooth function whose domain of definition contains $\bar A$. Analogously,
$B^-_n: \{z=z^-_n(v), u=0\}$
for a smooth function $z^-$. Thus, the points on $B^+_n$ and $B^-_n$ are uniquely determined by their $v$-coordinates.
This gives us a trivial projection of $B_n^+$ and $B_n^-$ to $A$, so we may formally consider the maps $F_n$, $n=0,\dots, N$,
as acting on $B_n^-$ or $B_n^+$ in the same way these maps act on $A$.

Since the foliations $E^{ss}$ and $E^{uu}$ are invariant with respect to the map $\Phi$, it follows that
$\Phi(E^{uu}_v)=E^{uu}_{F_0(v)}$ and $\Phi(E^{ss}_v)=E^{ss}_{F_0(v)}$
and, consequently,
the points $F_0^{-m_-}\circ\pi^{u}_n(x)$ and $\Phi^{-m_-}(x)$ have the same $v$-coordinate given any $x\in B_n$, and the same
is true for the points $F_0^{m_+}\circ\pi^{s}_n(x)$ and $\Phi^{m_+}(x)$. Thus, in the $v$-coordinates, we have
\begin{equation}\label{glcsc}
\Phi^{m_++m_-}|_{B^-_n} =F_0^{m_+}\circ F_n\circ F_0^{m_-}
\end{equation}

Denote as $T_n:(u,v,z)\mapsto(\bar u,\bar v,\bar z)$ the map $\Phi^{m_++m_-}$ from a small neighbourhood
of $B^-_n$ to a small neighbourhood of $B^+_n$.
The transversality condition implies that the image by the map $T_n$ of any leaf $(u=0,v=const)$
(of the foliation $E^{uu}$ in $W^u_{loc}$) is transverse to
$W^s_{loc}:\{\bar z=0\}$, which means that the derivative $\partial\bar z/\partial z$ is invertible.
Therefore, given any small $(u,\bar z)$ and $v\in \bar A$ we have a uniquely defined $(\bar u, z, \bar v)$ such that
$(\bar u,\bar z,\bar v)=T_n(u,z,v)$. So, we may write the map $T_n$ in the following form:
\begin{equation}\label{t1m0}
\bar u=p_n(u,v,\bar z), \quad \bar v = G_n(u,v,\bar z)=\bar F_n(v) + f_n(u,v,\bar z), \quad z=q_n(u,v,\bar z),
\end{equation}
where $p_n$, $q_n$, $f_n$ are smooth functions defined for small $(u,\bar z)$ and for $v$ from a small
neighbourhood $A''$ of $\bar A$ in $A$. We assume
\begin{equation}\label{faud}
f_n(0,v,0)\equiv 0
\end{equation}
in (\ref{t1m0}). As $u=0$ corresponds to the initial point in $W^u_{loc}$, and $\bar z=0$ corresponds
to the image of this point (by $T_n$) that lies in $W^s_{loc}$, having both $u=0$ and $\bar z=0$ corresponds to the initial
point lying in $B^-_n$ and its image lying in $B^+_n$. Thus, by (\ref{glcsc}), we have
\begin{equation}\label{nscmp}
\bar F_n=\Phi^{m_++m_-}|_{B^-_n} =F_0^{m_+}\circ F_n\circ F_0^{m_-},
\end{equation}
where $F_n$ is the scatering map. Since the cylinder $\bar A$ is invariant with respect to $F_0$, the maps $\bar F_n$ are defined
in a neighbourhood of $\bar A$, as the scattering maps $F_n$ are. Thus, we will further assume
that the open neighbourhood $A''$ of $\bar A$ in $A$ is chosen
such that the {\em modified scattering maps} $\bar F_n$
are all defined there and are homotopic to identity diffeomorphisms of $A''$ into $A$, moreover
\begin{equation}\label{appapap}
F_0^{-m_-}(A')\subseteq A''
\end{equation}
where $A'$ is a small neighbourhood of $\bar A$ in $A$ where the scattering maps $F_n$ are defined.

We recall that we assume
that the invariant cylinder $A$ persists for the class of perturbations of the map $\Phi$ we consider here (the perturbations are the maps $\tilde\Phi$
from a small neighbourhood $\mathcal V$ of $\Phi$ in the set of real-analytic exact symplectic maps), and $A$ depends continuously on the map.
Its stable and unstable manifolds, as well as the corresponding strong-stable and strong-unstable foliations, also depend on the map $\tilde\Phi$
continuously (in $C^1$). The transversality condition [S1] implies that the $C^1$-smooth homoclinic cylinders $B_1,\dots, B_N$ also persist
and depend continuously on the map. Thus, for all maps from $\mathcal V$, the scattering maps $F_1,\dots, F_N$ are defined
on the same open subset $A'$ of $A$ and the modified scattering maps $\bar F_1,\dots, \bar F_N$ are defined
on an open subset $A''$ of $A$ such that (\ref{appapap}) holds. All these maps are exact symplectomorphisms, homotopic to identity.
The maps $T_n$, $n=1,\dots, N$, remain well-defined for all $\tilde\Phi$ from $\mathcal V$, and formulas (\ref{t1m0})-(\ref{nscmp}) hold.
Estimates (\ref{mpc})-(\ref{derv}) also remain valid (as they are true for every map with a compact,
symmetrically normally-hyperbolic invariant cylinder $A$). In what follows we will consider arbitrary maps from $\mathcal V$ (e.g.
we do not assume the invariance of the cylinder $\bar A$ for all such maps; from now on it is just a cylinder bounded by two arbitrary,
non-intersecting, essential curves $\gamma^\pm\subset A'$).

\medskip

Let us denote as $T_0$ the map $\tilde\Phi$ restricted to the $\delta$-neighbourhood of $A$. Let us call
the union of the $\delta$-neighbourhood
of $A$ with certain small neighbourhoods of the cylinders $\tilde\Phi (B^-),\dots, \tilde\Phi^{m_++m_--1}(B^-)$
{\em a homoclinic channel}. For every finite orbit in the homoclinic channel there is a uniquely defined sequence of points
$P_s$ ($s=0,\dots, 2J+1$) of this orbit which lie in the $\delta$-neighbourhood of $A$ and satisfy
$P_{2j+1}=T_0^{k_j}P_{2j}$ ($j=0,\dots,J$) and $P_{2j}=T_{n_j}P_{2j-1}$ ($j=1,\dots,J$), where $n_j$
may take values from $1,\dots,N$ and $k_j\geq 0$. We will call the sequence $P_j$ a {\em channel orbit}, and the sequence
$k_0,\{n_s,k_s\}_{1\leq s \leq J}$ will be called {\em the code} of the orbit.
Given a code $k_0,\{n_s,k_s\}_{1\leq s \leq J}$, we may consider a {\em shadow orbit}, namely a sequence
$v_s^*$ of points in $A$ such that $v_{2j+1}^*=F_{0}^{k_j}v_{2j}^*$ and $v_{2j}^*=\bar F_{n_j}v_{2j-1}^*$. Note that
we always assume that
\begin{equation}\label{dfshdd}
v^*_{2j-1}\in A'', \qquad j=1,\dots, J,
\end{equation}
in order to have the maps $\bar F_{n_j}$ well-defined (this, in particular, implies that not all codes necessarily correspond to a shadow orbit).
Given a channel orbit $P_s$, its shadow is the shadow orbit of the corresponding code
with $v_0^*$ equal to the $v$-coordinate of $P_0$.

\subsection{Shadowing orbits of proper codes}

Our next goal is to estimate the deviation of the channel orbit $P_s$ from its shadow. Here, we restrict our attention
onto a special class of orbits which correspond to a special class of codes.
Namely, we will call a finite code {\em proper} if the corresponding sequence $k_s$ satisfies for every $s$
\begin{equation}\label{oper}
k_s\geq \bar k \mbox{~~and~~~} k_s\geq \gamma k_{s+1} + D,
\end{equation}
for some $\bar k \geq 0$, $D\geq 0$ and $\gamma>1$. In other words, $k_s$ is a sufficiently fast decreasing sequence
of sufficiently large numbers.

\begin{lemma}\label{shadprop}
Given any sufficiently large $\bar k$, $\gamma$ and $D$, for any shadow orbit $v_0^*, \dots, v^*_{_{2J+1}}$
with a proper code $k_0,\{n_s,k_s\}_{1\leq s \leq J}$, given any
$u^{in}$ and $z^{out}$ such that $\|u^{in}\|\leq\delta$, $\|z^{out}\|\leq\delta$, in the $\delta$-neighbourhood of $A'$
there exists a uniquely defined channel orbit $P_s(u_s,v_s,z_s)$ such that $u_0=u^{in}$, $v_0=v_0^*$,
$z_{_{2J+1}}=z^{out}$, and $P_{2j+1}=T_0^{k_j}P_{2j}$,
$P_{2j}=T_{n_j}P_{2j-1}$. Moreover,
\begin{equation}\label{shade}
\|v_s-v_s^*\|\leq \delta(\alpha\lambda)^{k_{_J}/2}\leq 2\delta(\alpha\lambda)^{\bar k/2},
\end{equation}
and
\begin{equation}\label{shind}
\|u_{_{2J+1}}\|\leq \delta \lambda^{\bar k}, \qquad \|z_0\|\leq \delta \lambda^{\bar k}.
\end{equation}
\end{lemma}
\begin{remark} {\em Usual shadowing results would require hyperbolicity (or its topological analogues) from
the maps $F_0$ and $F_1,\dots,F_N$, see e.g. \cite{DG2013}. We, however, do not make any assumption on the dynamics of these maps in this
lemma (e.g. we have not assumed the symplecticity so far). Therefore we need to restrict here the class of shadow orbits
to those with proper codes only;
we believe any significantly stronger shadowing statement can not hold in this situation without further assumptions.}
\end{remark}

\noindent
{\em Proof of the lemma}.
For $J=0$ the statement of the lemma is contained in Lemma~\ref{Le:1},
so we will proceed by induction in $J$. Thus, we may make an assumption that
given any $\tilde z$ such that $\|\tilde z\|\leq \delta$ the sequence
$(u_s,v_s,z_s)$ with the code $k_0,n_1,\dots,k_{J-1}$ at $s=0,\dots,2J-1$ is uniquely defined by the condition
$u_0=u^{in}$, $v_0=v_0^*$ and $z_{_{2J-1}}=\tilde z$, and that the condition
\begin{equation}\label{indshade}
\|v_s-v_s^*\|\leq \delta(\alpha\lambda)^{k_{_{J-1}}/2}
\end{equation}
is fulfilled for all $s\leq 2J-1$. Denote $u_{_{2J-1}}=\tau(\tilde z)$ and
$v_{_{2J-1}}=\phi(\tilde z)$ (we assume $u^{in}$ and $v_0^*$ fixed, and do not indicate the dependence on them).
Note that it follows from (\ref{dfshdd}),(\ref{indshade}) that
\begin{equation}\label{arrho}
\phi(\tilde z)\in A''_\rho \mbox{~~~for any~~} \rho>\delta(\alpha\lambda)^{\bar k/2},
\end{equation}
where $A''_\rho$ is the closed $\rho$-neighbourhood of $A''$.

Since $(u_{_{2J-1}},v_{_{2J-1}},\tilde z)=T_0^{k_{_{J-1}}}(u_{_{2J-2}},v_{_{2J-2}},z_{_{2J-2}})$, it follows from
Lemma~\ref{Le:1} (see (\ref{lamb})) that
\begin{equation}\label{taud}
\|\tau\|\leq\delta\lambda^{k_{_{J-1}}}\leq \delta\lambda^{\bar k}.
\end{equation}
We will also include in our induction assumption a bound for the derivatives:
\begin{equation}\label{tdr}
\|\tau',\phi'\|\leq\nu
\end{equation}
for some sufficiently small constant $\nu$. Thus, in order to carry out the
induction, when we prove that the sought sequence $(u_j,v_j,z_j)$ is uniquely defined for all $j=0,\dots,2J+1$
we must also show that
\begin{equation}\label{tdrz}
\|\partial (u_{_{2J+1}},v_{_{2J+1}})/\partial z^{out}\|\leq\nu
\end{equation}
with the same $\nu$.

Since $(u_{_{2J+1}},v_{_{2J+1}},z^{out})=T_0^{k_{_J}}(u_{_{2J}},v_{_{2J}},z_{_{2J}})$, it follows from Lemma~\ref{Le:1} that
$z_{2J}$ is a uniquely defined smooth function of $(u_{_{2J}},v_{_{2J}})$ and $z^{out}$; we denote
this function as $\sigma: (u_{_{2J}},v_{_{2J}},z^{out})\mapsto z_{_{2J}}$. By (\ref{lamb}),(\ref{dery}),
\begin{equation}\label{saud}
\|\sigma\|\leq\delta\lambda^{k_{_J}}\leq\delta\lambda^{\bar k},
\end{equation}
and
\begin{equation}\label{sdrt}
\|\sigma'\|\leq\nu,
\end{equation}
for any $\nu$ chosen in advance (if $\bar k$ is large enough), and
\begin{equation}\label{zmso}
\|\partial\sigma/\partial z^{out}\|\leq \lambda^{k_{_J}}.
\end{equation}

With these notations, the sought sequence $(u_j,v_j,z_j)$ exists indeed and is defined uniquely if
and only if the following equation has a unique solution
\begin{equation}\label{lsts}
\begin{array}{l}\displaystyle
u_{_{2J}}=p_{n_{_J}}(\tau(z_{_{2J-1}}),\phi(z_{_{2J-1}}),\sigma (u_{_{2J}},v_{_{2J}},z^{out})),\\ \\ \displaystyle
v_{_{2J}}=G_{n_{_J}}(\tau(z_{_{2J-1}}),\phi(z_{_{2J-1}}),\sigma (u_{_{2J}},v_{_{2J}},z^{out})),\\ \\\displaystyle
z_{_{2J-1}}=q_{n_{_J}}(\tau(z_{_{2J-1}}),\phi(z_{_{2J-1}}),\sigma (u_{_{2J}},v_{_{2J}},z^{out})),
\end{array}
\end{equation}
(see (\ref{t1m0})). By construction, the functions $p_n$, $q_n$, $G_n$ in (\ref{t1m0}) are defined
at all $u$ and $\bar z$ which are sufficiently small and all $v\in A''$. Obviously, they will remain defined and
smooth for $v$ from a sufficiently small neighbourhood of $A''$. Thus, to be sure that the system (\ref{lsts})
is well-defined, we must check that, by taking $\bar k$ sufficiently large, the values of $\tau$ and $\sigma$ can be made arbitrarily small
and the range of values of $\phi$ can be confined to an arbitrarily small neighbourhood of $A''$,
and this is indeed given by (\ref{arrho}),(\ref{taud}),(\ref{saud}).

As the functions $p_n$, $q_n$, $G_n$ are smooth, their derivatives are
bounded by a constant. We, therefore, write
$$\|p_n',q_n',G_n'\|\leq C$$
and take $\nu$ in (\ref{tdr}),(\ref{sdrt}) such that
$$C\nu<1.$$
Thus, the right-hand side of system (\ref{lsts}) is a contracting operator, which will immediately give us
the existence and uniqueness of the solution once we check that given any
$z$, $u$, $v$ and $z^{out}$ such that $\|z\|\leq\delta$, $\|z^{out}\|\leq\delta$, $\|u\|\leq\delta$
and $v\in A''_\rho$ (for some $\rho$ small enough) the functions $p_n(\tau(z),\phi(z),\sigma (u,v,z^{out}))$, $G_n(\tau(z),\phi(z),\sigma (u,v,z^{out}))$,
$q_n(\tau(z),\phi(z),\sigma (u,v,z^{out}))$ return the new values of, respectively, $u$, $v$ and $z$ for which
the same conditions are satisfied (namely, $\|z\|\leq\delta$, $\|u\|\leq\delta$
and $v\in A''_\rho$). That $\|p_n\|\leq\delta$ and $\|q_n\|\leq\delta$, this follows immediately, as
$p_n$ and $q_n$ define the map $T_n$ which acts from a small neighbourhood of the cylinder $B^-_n$ to a
small neighbourhood of the cylinder $B^+_n$, and both the cylinders lie inside the $\delta$-neighbourhood of $A$.
So, we are left to verify that $v_{_{2J}}=G_{n_J}(\tau,\phi,\sigma)\in A''_\rho$.
Note that by the induction assumption (\ref{indshade})
we have $\|\phi-v^*_{_{2J-1}}\|\leq 2\delta(\alpha\lambda)^{k_{_{J-1}}/2}$. Since $G_n:=\bar F_n+f_n$, and
$v^*_{_{2J}}=\bar F_{n_{_J}}(v^*_{_{2J-1}})$, it follows that
$$\|G_{n_{_J}}(\tau,\phi,\sigma)-v^*_{_{2J}}\|\leq C \|\phi-v^*_{_{2J-1}}\|+\|f_{n_{_J}}(\tau,\phi,\sigma)\|,$$
hence
\begin{equation}\label{vmeq}
\|v_{_{2J}}-v^*_{_{2J}}\|\leq C\delta (2(\alpha\lambda)^{k_{_{J-1}}/2}+\lambda^{k_{_J}})
\end{equation}
(we take into account that $f_n$ vanishes, by (\ref{faud}), at $\tau=0$, $\sigma=0$, so $\|f_n\|\leq C\|\tau,\sigma\|
\leq C\lambda^{k_{_J}}$ by (\ref{taud}),(\ref{saud})). At $\bar k$ large enough, we have from (\ref{vmeq}) that
$$\|v_{2J}-v^*_{2J}\|\leq \rho.$$
Since $v^*_{2J}$ lies in $A''$, it follows that
$v_{_{2J}}=G_{n_{_J}}(\tau,\phi,\sigma)\in A''_\rho$ indeed.

Thus, system (\ref{lsts}) has a unique solution $(u_{_{2J}},v_{_{2J}},z_{_{2J-1}})$ as required; moreover,
by (\ref{zmso}) we have
\begin{equation}\label{thisp}
\|\partial (u_{_{2J}},v_{_{2J}})/\partial z^{out}\|\leq \frac{C}{1-C\nu}\lambda^{k_{_J}}=o(\alpha^{-k_{_J}})
\end{equation}
(as $\alpha\lambda<1$ by (\ref{alla})). Since $(u_{_{2J+1}},v_{_{2J+1}},z^{out})=
T_0^{k_{_J}}(u_{_{2J}},v_{_{2J}},z_{_{2J}})$,
we find from (\ref{dery}),(\ref{derz}) and(\ref{thisp}) that
$\|\partial (u_{_{2J+1}},v_{_{2J+1}})/\partial z^{out}\|$ can be made as small as we need by taking $k_{_J}$ large enough,
i.e. (\ref{tdrz}) holds true indeed for $\bar k$ large enough.

We have proved the existence and uniqueness of the sought sequence
$(u_s,v_s,z_s)$, $s=0,\dots,2J+1$.
It remains to prove inequality (\ref{shade}). At $j\leq 2J-1$ this is given by the induction assumption
(\ref{indshade}), since $k_{_{J-1}}>k_{_J}$. In order to check (\ref{shade}) at $j\geq 2J$, we recall that
$(u_{_{2J+1}},v_{_{2J+1}},z_{_{2J+1}})=T_0^{k_{_J}}(u_{_{2J}},v_{_{2J}},z_{_{2J}})$ and
$v^*_{_{2J+1}}=F_0^{k_{_J}} v^*_{_{2J}}$. It then follows from (\ref{lf00}), (\ref{vmeq}) and (\ref{shad0}) that
$$\|v^*_{_{2J+1}}-F_0^{k_{_J}} v_{_{2J}}\|\leq C\delta (2(\alpha\lambda)^{k_{_{J-1}}/2}+\lambda^{k_{_J}})
\alpha^{k_{_J}},$$
and
$$\|v_{_{2J+1}}-F_0^{k_{_J}} v_{_{2J}}\|\leq \delta (\alpha\lambda)^{k_{_J}/2},$$
which gives
$$\|v_{_{2J+1}}-v^*_{_{2J+1}}\|\leq \delta (\alpha\lambda)^{k_{_J}/2}(1+C(\alpha\lambda)^{\bar k/2})+
2C\delta (\alpha\lambda)^{k_{_{J-1}}/2}\alpha^{k_{_J}}.$$
We finish the proof of the lemma by noting that this inequality and inequality (\ref{vmeq}) imply (\ref{shade}) at $j=2J+1$ and, respectively, $j=2J$ if
$$\gamma > \ln\frac{1}{\alpha\lambda}\left/\ln\frac{\alpha}{\lambda}\right., \qquad
D> 2\ln \frac{2C}{1-C(\alpha\lambda)^{\bar k/2}}\left/\ln\frac{\alpha}{\lambda}\right.$$
in (\ref{oper}). Inequality (\ref{shind}) follows immediately from (\ref{lamb}). ~~$\sqcap\!\!\!\!\sqcup$

\subsection{Replacing a code with a proper code}\label{s43}

If the diffeomorphism $F_0$ is area-preserving, the Poincare Recurrence Theorem implies that recurrent (Poisson stable) orbits
of $F_0$ are dense in the invariant cylinder $A$. This fact, as the following lemma shows,
allows an arbitrary orbit of the iterated function system $\{F_0, F_1, \dots, F_N\}$ to be approximated by a shadow with a proper code.

\begin{lemma}\label{pncccl}
Let $v_0,\dots, v_{_{2J+1}}$ be a sequence of points in $int(A)$ such that
\begin{equation}\label{shadorb}
v_{2j+1}=F_0^{i_j}(v_{2j}) \;\; (j=0,\dots,J), \qquad v_{2j}=F_{n_j}(v_{2j-1})\;\; (j=1,\dots,J)
\end{equation}
for some $i_j\geq 0$ and $n_j=1,\dots,N$ (we assume that the points $v_{2j-1}$, $j=1,\dots,J$, belong to the open set $A'$,
the domain of definition of the maps $F_n$ with $n\geq 1$).
Then, given two open subsets of $A$, an arbitrary neighbourhood $U_0$ of the point $v_0$ and an arbitrary neighbourhood $U_{_{2J+1}}$ of the point $v_{_{2J+1}}$,
and given any $\bar k$, $\gamma$ and $D$, in $int(A)$ there exists a sequence of points $v_s^*$ such that
$$v_0^*\in U_0, \quad v^*_{_{2J+1}}\in U_{_{2J+1}}, \qquad
v^*_{2j+1}=F_0^{k_j}(v^*_{2j}) \;\; (j=0,\dots,J), \qquad v^*_{2j}=\bar F_{n_j}(v^*_{2j-1})\;\; (j=1,\dots,J)$$
with the same $n_j$ as in (\ref{shadorb}), $v^*_{2j-1}\in A''$ (the domain of definition of the maps $\bar F_n$) for $j=1,\dots,J$,
and the numbers $k_j$ form a proper sequence in the sense of (\ref{oper}).
\end{lemma}

\noindent
{\em Proof.} By the definition of the modified scattering maps $\bar F_n$ (see (\ref{nscmp})), it is enough to show that in $int(A)$
there exists a sequence of points $\hat v_s$
such that $\hat v_0\in U_0$, $\hat v_{_{2J+1}}\in U_{_{2J+1}}$,
\begin{equation}\label{vdj}
\hat v_{2j+1}=F_0^{\hat k_j}(\hat v_{2j}) \;\; (j=0,\dots,J), \qquad \hat v_{2j}=F_{n_j}(\hat v_{2j-1})\;\; (j=1,\dots,J)
\end{equation}
with the same $n_j$ as in (\ref{shadorb}), $\hat v_{2j-1}\in A'$ for $j=1,\dots,J$, and the numbers $\hat k_j$ are such that
$\hat k_J=k_J+m_+$, and, if $J\geq 1$, then $\hat k_j=k_j+m_++m_-$ at $0\leq j \leq J-1$,
where $k_j$, $j=0,\dots,J$, form a proper sequence. The sought sequence $v^*_s$ will then be given by
$$v_0^*=\hat v_0, \qquad v_{_{2J+1}}^*=\hat v_{_{2J+1}}, \qquad v^*_{2j-1}=F_0^{-m_-} \hat v_{2j-1}, \qquad v^*_{2j}=F_0^{m_+} \hat v_{2j} \quad (j=1,\dots, J);$$
we will have $v_{2j-1}^*\in A''$ by virtue of (\ref{appapap}).

We will proceed by induction in $J$. Let $\hat U$ be a small neighbourhood of $v_1$ in $A$ such that $F_0^{-i_0}\hat U\subset U_0$. If $J=0$ we assume $\hat U\subseteq U_1$;
if $J\geq 1$, we assume that $\hat U$ is small enough,
so that the map $F_{n_1}$ is defined everywhere on $\hat U$. By induction, there is a point $v^\prime\in U_2=F_{n_1}\hat U$ such that
${\mathcal F} v^{\prime} \in U_{_{2J+1}}$, where
$\displaystyle \mathcal F=\left (\prod_{2\leq j\leq J} F_0^{\hat k_j}\circ F_{n_j}\right) \circ F_0^{\hat k_1}$
with the numbers $\hat k_j$ such that, if we define
\begin{equation}\label{dfht}
k_J=\hat k_J-m_+ \mbox{~~~and,~~at~~~} J\geq 2\mbox{,~~~~~} k_j=\hat k_j-m_--m_+ \mbox{~~~~for all~~~} j\leq J-1,
\end{equation}
then the numbers $k_j$, $j=1,\dots,J$, form a proper sequence. Let $U'\subseteq \hat U$ be a small enough neighbourhood of $v_1$ in $A$ such that
$\mathcal F\circ F_{n_1} (U')\subseteq U_{_{2J+1}}$; at $J=0$ take $U'=\hat U\subseteq U_1$.

By construction, we will prove the lemma if we find a point $\hat v_1\in U'$ such that $\hat v_0=F_0^{-\hat k_0} \hat v_1\in U_0$ for some $\hat k_0$
such that
\begin{equation}\label{reqsa}
\hat k_0-m_--m_+\geq \bar k
\end{equation}
and, at $J\geq 1$,
\begin{equation}\label{reqsb}
\hat k_0-m_--m_+\geq \gamma k_1+D,
\end{equation}
with $k_1$ given by (\ref{dfht}). Indeed, if we define $k_0=\hat k_0-m_+$ at $J=0$ and $k_0=\hat k_0-m_+-m_-$ at $J\geq 1$, then
the sequence $k_0,\dots,k_J$ will be proper (see (\ref{oper})), and for the corresponding sequence $\hat v_s$ defined by (\ref{vdj})
we will have $\hat v_{_{2J+1}}\in U_{_{2J+1}}$, as $\hat v_{_{2J+1}}=\mathcal F\circ F_{n_1} (\hat v_1)\in \mathcal F\circ F_{n_1} (U')\subseteq U_{_{2J+1}}$ at $J\geq 1$
and $\hat v_{_{2J+1}}=\hat v_1\in U'\subseteq U_1$ at $J=0$.

Now note that since $F_0$ is an area-preserving map of a bounded region $A$, recurrent orbits are dense in $A$. In particular, there is a recurrent point in $U'$.
It is easy to see that any such point can serve as the sought point $\hat v_1$. Indeed, if $\hat v_1\in U'$ is recurrent,
then one can find an arbitrarily large $k$ such that $F_0^{-k}\hat v_1\in U'$. As $U'\subseteq \hat U$, it follows that $F_0^{-i_0-k} \hat v_1\in U_0$.
So, we put $\hat k_0=k+i_0$, and conditions (\ref{reqsa}),(\ref{reqsb}) are fulfilled if the return time $k$ was taken large enough.
\hfill$\Box$

We say that two points $v_0$ and $v_m$ are connected by an orbit of the iterated function system $\{F_0,\ldots,F_N\}$ if $v_{_{2J+1}}$ is an image of $v_0$
by a certain sequence of maps $F_n$. Obviously, this means that $v_0$ and $v_m$ are the first and the last points in a sequence of points $v_s$ constructed
by the rule (\ref{shadorb}) with $m=2J+1$. Since the corresponding sequence $v^*_s$ constructed in Lemma \ref{pncccl} is a shadow of proper code, we may apply Lemma \ref{shadprop}
to it. Thus, combining Lemmas \ref{pncccl} and \ref{shadprop}, we obtain the following

\begin{lemma}\label{lmsw}
Let the map $F_0$ be area-preserving. Let two points $v_0\in A$ and $v_m\in A$ be connected by an orbit of
the iterated function system $\{F_0,\ldots,F_N\}$.
The, for any $\varepsilon>0$ the $\varepsilon$-neighbourhoods
of $v_0$ and $v_m$ in $\mathbb R^{2d}$ are connected by an orbit of the map $\tilde\Phi$.
\end{lemma}

\section{Transport in the homoclinic channel in the symplectic case}

\subsection{Symplecticity of invariant manifolds and scattering maps}\label{s51}

In this Section we continue the analysis of the behaviour in the homoclinic channel. Now we take into account the fact
that the map $\tilde\Phi$ under consideration is exact symplectic.
We start with establishing some useful geometric properties of the stable and unstable
manifolds and the scattering maps. These properties are based on a symplectic orthogonality property, as given by the
following

\begin{proposition}\label{Pro:ort}
If $A$ is a symmetrically normally-hyperbolic invariant manifold and $x\in A$, then
$T_yW^{s}(A)\perp_{\Omega}T_yE^{ss}(x)$
for any $y\in E^{ss}(x)$ and $T_yE^{uu}(x)\perp_{\Omega}T_yW^u(A)$ for any $y\in E^{uu}(x)$.
\end{proposition}

\noindent{\em Proof.}
Let $y\in W^s(A)$. Take any $w\in T_yE^{ss}(x)$ and $u\in T_yW^s(A)$.
Since the map $\tilde\Phi$ preserves the form $\Omega$, we have for any $m\in\mathbb N$:
$$\Omega(w,u)= \Omega((\tilde\Phi')^m w, (\tilde\Phi')^m u)=(\alpha\lambda)^m \Omega(\alpha^{-m}(\tilde\Phi')^m w, \lambda^{-m}(\tilde\Phi')^m u)=
O((\alpha\lambda)^m)\,.
$$
Taking the limit $m\to+\infty$, we find that $\Omega(w,u)=0$, i.e. $u\perp_{\Omega} v$. Thus, we have proved
$T_yE^{ss}(x)\perp_{\Omega}T_yW^s(A)$.
In a similar way we conclude that $T_yE^{uu}(x)\perp_{\Omega}T_yW^u(A)$ for any $y\in W^u(A)$.
\hfill$\Box$

Note that we do not use that the manifold $A$ is a two-dimensional cylinder embedded
in $\mathbb R^{2d}$. In fact, this and other statements we prove here (Propositions \ref{smplca}-\ref{Prop:sympl})
assume only that $A$ is a submanifold of a smooth manifold $M$ endowed with a closed non-degenerate symplectic form $\Omega$,
and that the diffeomorphism $\tilde\Phi$ preserves $\Omega$.

\medskip

\begin{proposition}\label{smplca}
The restriction of the symplectic form
$\Omega$ to the symmetrically normally-hyperbolic invariant manifold $A$ is non-degenerate.
\end{proposition}
\noindent
{\em Proof.}
If the proposition is not true and the restriction of the symplectic
form is degenerate, then there is a non-zero vector $w\in T_x(A)$
such that $w\perp_{\Omega} T_x(A)$. On the over hand
$w\in T_xA=T_xW^s(A)\cap T_xW^u(A)$ implies that
$w\perp_{\Omega} T_x E^{ss}_x$ and $w\perp_{\Omega} T_x E^{uu}_x$.
The normal hyperbolicity assumptions imply that
$TM=T_xE^{ss}_x\oplus T_xE^{uu}_x\oplus T_xA$ for any $x\in A$.
Consequently, $w\perp_{\Omega} T_xM$, which contradicts
to the non-degeneracy of $\Omega$, and the proposition follows immediately.
\hfill$\Box$

\bigskip

We call a homoclinic intersection of $W^u(A)$ and $W^s(A)$ at a point $y$ strongly transverse
if $E^{uu}_y$ is transverse to $W^s(A)$ and $E^{ss}_y$ is transverse to $W^u(A)$ at the point $y$.

\begin{proposition}
If $y\in E^{uu}(x_1)\cap E^{ss}(x_2)$ for some $x_1,x_2\in A$ and
$TM=T_yE^{uu}(x_1)\oplus T_yW^s(A)$ then $TM=T_yE^{ss}(x_2)\oplus T_yW^u(A)$
and, consequently, the homoclinic intersection at $y$ is strongly transverse.
\end{proposition}

The proof of this proposition is completely straightforward: it is sufficient to note that any vector from $T_yE^{ss}(x_2)\cap T_yW^u(A)$
is $\Omega$-orthogonal to all vectors due to Proposition~\ref{Pro:ort}.
The proposition implies that the strong transversality is equivalent to
the transversality of the strong stable leaves to the unstable manifold (or the transversality of the strong unstable leaves
to the stable manifold). This property reduces the number of conditions which are necessary to verify the strong transversality
of a homoclinic intersection.

\medskip

For every $y\in W^s(A)$ there is a unique $x\in A$ such that $y\in E^{ss}(x)$. We define the projection $\pi^s:W^s(A)\to A$
by setting $\pi^s(y)=x$. Let $u_1,u_2,\dots$ be some coordinates on $A$ defined in a small neighbourhood $U$ of the point $x$.
Define coordinates $(u,v)$ in $(\pi^{s})^{-1}(U)$ such that $v=0$ corresponds to a point in $A$ and
$u=\mathrm{const}$ corresponds to the strong stable leave $E^{ss}_{(p,q,0)}$.
In these coordinates $\pi^s:(u,v)\mapsto (u,0)$.

Since $T_yE^{ss}(x)\perp_{\Omega}T_yW^s(A)$, we see that
in these coordinates $\Omega|_{W^s(A)}=\sum_{i,j} a_{ij}(u,v)du_i\wedge du_j$. On the other hand, the symplectic form is closed, i.e., $d\Omega=0$.
We have $d\Omega|_{W^s(A)}=\sum_{i,j}\frac{\partial a_{ij}}{\partial v}dv\wedge du_i\wedge du_j$, which implies that the coefficients
$a_{ij}$ do not depend on $v$, and consequently $\Omega|_{W^s(A)}=\sum_{i,j} a_{ij}(u)du_i\wedge du_j$.

Let $B$ be any section of $W^s(A)$ transverse to the strongly stable leaves. Then
the restriction $\pi^{s}|_{B}:B\to A$  is a local diffeomorphism. Moreover,
since the projection is the identity in the coordinates $u$, we find that
$\pi^s|_{B}$ is a symplectomorphism, i.e. it transforms $\Omega|_{B}$ into $\Omega|_A$. In particular, $\Omega|_{B}$
is non-degenerate, i.e. $B$ is a symplectic manifold.

Obviously, a similar statement is true for the stable manifolds replaced by the unstable ones:
for any section $B$ of $W^u(A)$ transverse to the strongly unstable leaves, the projection
$\pi^u: B\to A$ by the strongly unstable leaves is locally a symplectomorphism.
Thus, we obtain the following

\begin{proposition}
If $y\in W^s(A)\cap W^u(A)$ is a strongly transverse homoclinic point and
$B$ is a sufficiently small neighbourhood of $y$ inside $W^s(A)\cap W^u(A)$, then
the scattering map $F_{B}=\pi^s|_{B}\circ (\pi^u|_{B})^{-1}:B^u\to B^s$ is a symplectomorphism,
where $B^{u,s}=\pi^{u,s}(B)\subset A$.
\end{proposition}

We can define the scattering map $F_B$ relative to any connected subset $B$ of $W^s(A)\cap W^u(A)$ that consists of strongly transverse homoclinic points.
When $B$ is not a small neighbourhood of a single point, the scattering map $F_B$ does not need to be single-valued nor injective
(eventhough every branch of it is a local diffeomorphism). In the example of simple homoclinic cylinder we consider in this paper,
the scattering map is, however, single-valued and injective (so it is a symplectic diffeomorphism for a large open subset $A'$ of $A$).

Assume the symplectic form is exact: $\Omega=d\vartheta$, where $\vartheta$ is a differential 1-form (in the
case of our interest, $M=\mathbb R^{2d}$, $\Omega=dp\wedge dq$, and $\vartheta=p dq$). The symplectic map $\tilde\Phi$ is
exact when
$$\int_\gamma \vartheta=\int_{\tilde\Phi(\gamma)}\vartheta$$
for every smooth closed curve $\gamma$. Obviously, the exactness of $\tilde\Phi$ implies the exactness of the map $F_0=\tilde\Phi|_A$.

\begin{proposition}\label{Prop:sympl}
Let $A'\subseteq A$ be a region such that the scattering map
$F_B$ is a diffeomorphism $A'\to F_B(A')$; moreover, for each point $x\in A'$ the corresponding leaves $E^{uu}(x)$ and $E^{ss}(F_B(x))$ intersect
$B$ exactly at one point. Then the restriction of $F_B$ on $A'$ is an exact map.
\end{proposition}
\noindent
{\em Proof.} Let us prove that the map $(\pi^u|_B)^{-1}$ is exact on $A'$. The proof of the exactness of the map $(\pi^s|_B)^{-1}$ on $F_B(A')$
is exactly the same, so the exactness of $F_B$ will follow immediately.
Take any smooth closed curve $\gamma\subset A'$. By assumption, for any $x\in\gamma$ there is a unique point $y(x)\in B$ such that $y\in E^{uu}(x)$,
the union of the points $y(x)$ over all $x\in\gamma$ gives the curve $(\pi^u|_B)^{-1}\gamma=\tilde\gamma\subset B$. As the strongly unstable leaves
are simply-connected (each is a diffeomorphic copy of $\mathbb R^k$ where $2k=dim (M)-dim(A)$) and depend smoothly on the base point $x$, one can connect
each point $x\in\gamma$ with the corresponding point $y(x)\in\tilde\gamma$ by a smooth arc $\ell(x)$ that lies in $E^{uu}(x)$ so that the union of these arcs forms
a smooth two-dimensional surface $S\subset W^u(A)$, an annulus bounded by $\gamma$ and $\tilde\gamma$. By Stokes theorem,
$$\int_\gamma \vartheta-\int_{\tilde\gamma}\vartheta=\int_S \Omega.$$
At every point $y\in S$ the tangent plane contains the vector tangent to one of the curves $\ell(x)$ which lies in the $E^{uu}(x)$, so $\Omega$ vanishes on $T_yS$ by Proposition
\ref{Pro:ort}. Thus, $\int_S \Omega=0$, which gives us the required identity $\int_\gamma \vartheta=\int_{\tilde\gamma}\vartheta$ for every smooth closed curve
$\gamma$ in $A'$.
\hfill$\Box$

Note that, surprisingly, the exactness of the scattering map in the statement above does not require the exactness of the map $\tilde\Phi$ itself.

\subsection{Transport in an iterated functions system and obstruction curves}

The symplecticity of the map $F_0=\tilde\Phi|_A$ established in Proposition \ref{smplca}
means that this map is area-preserving (with the area of a domain obtained by integrating $\Omega|_{A}$ over this domain).
Therefore, as shown in Section \ref{s43}, for two open sets to be connected by an orbit from the homoclinic channel
it is enough for these sets to be connected by the orbits of the iterated function system $\{F_0, F_1, \dots, F_N\}$. As we showed in Section \ref{s51}
all these maps are exact symplectomorphisms. The diffeomorphism $F_0$ is defined everywhere on the cylinder $A$ which is invariant with respect to $F_0$, i.e. $F_0(A)=A$.
The scattering maps $F_n$, $n=1,\dots,N$, are defined on a subset $A'$ of the cylinder $A$ and, as follows from the simplicity assumptions [S1]-[S3], they are homotopic
to identity diffeomorphisms $A'\to A$. The exact symplecticity of the maps $F_n$  implies that the area
between any curve $\gamma$ and its image $F_n(\gamma)$ is zero. Hence, $F_n(\gamma)\cap\gamma\ne\emptyset$
for any essential curve $\gamma\subset A'$.

We assume that there exist coordinates $v=(y,\varphi)$
in $A$ such that the map $F_0:(y,\varphi)\mapsto (\bar \varphi,\bar y)$
in these coordinates satisfies the {\em twist condition}, i.e.
$$\frac{\partial\bar \varphi}{\partial y}\neq 0$$
everywhere in this cylinder (we assume that $\varphi\in \mathbb S^1$ is the angular variable).
Recall that according to the Birkhoff theory, any invariant essential curve of a twist map is a graph of
a Lipschitz function $y(\varphi)$ \cite{Herman1983}.

Let $\bar A$ be a compact cylinder in $A'$ bounded by two simple essential curves $\gamma^+$ and $\gamma^-$ such that
$\gamma^-\cap\gamma^+=\emptyset$ (we no longer need to assume that $\bar A$ is an invariant cylinder). Let $\gamma^+$ corresponds
to larger values of $y$ than $\gamma^-$ does. The set $A\setminus int(\bar A)$ consists of two connected components, the upper component
$A^+$ contains $\gamma^+$ and the lower component $A^-$ contains $\gamma^-$. If $\bar A$ contains
an essential curve $\gamma^*$ which is {\em invariant under all of the
maps $F_n$}, $n=0, \dots, N$, then the curve $\gamma^*$ separates the cylinder $\bar A$ into two invariant parts, so no trajectory of the iterated
function system $\{F_0, F_1, \dots, F_N\}$ which starts in $A^-$ can get to $A^+$. In other words, the absence of essential common invariant curves
in $\bar A$ is a necessary condition for the orbits of iterated function system to connect $A^-$ with $A^+$. The following theorem shows that this
condition is also sufficient.

\begin{theorem}\label{Thm:Birkhoff}
Let $F_1,\ldots,F_{N}$ be exact symplectomorphisms $A'\to A$, homotopic to identity. Let $A$ be invariant with respect to a symplectic diffeomorphism
$F_0$ which satisfies the twist condition on $A$. Suppose no essential curve in $\bar A$ is invariant with respect to all maps $F_n$, $n=0,1,\dots, N$.
Then there is a finite trajectory  $(v_i)_{i=0}^m \subset \bar A$
of the iterated function system $\{F_0, F_1, \dots, F_N\}$ that starts on $\gamma^-$ and ends on $\gamma^+$ (i.e.
$v_0\in\gamma^-$, $v_m\in\gamma^+$, and and $v_{i+1}=F_{k_i}(v_i)$ for some sequence of $k_i=0,\dots,N$).
\end{theorem}

\begin{remark} {\em As the common invariant curve is, in particular, an invariant curve of the twist map $F_0$, the Birkhoff theory implies that it
is necessarily must be a graph of a Lipschitz function $y = y^*(\varphi)$, so it is sufficient to verify the absence of Lipschitz common invariant curves}.
\end{remark}

\begin{remark} {\em Our statement makes an important change
in the setup of the problem compared to e.g. \cite{M2002,L2007} as we do not ask the boundaries $\gamma^-$ and $\gamma^+$
to be invariant with respect to any of the maps $F_n$, $n=0,\dots,N$. Indeed, it is not natural to assume that the scattering maps
preserve boundaries as this would require certain non-transversality of stable and unstable manifolds
associated with $\tilde\Phi$-invariant curves on the boundary.}
\end{remark}

\medskip

\noindent
{\em Proof of Theorem \ref{Thm:Birkhoff}.}
The boundary of the $F_0$-invariant cylinder $A$ consists of two non-intersecting essential curves; we will refer to the boundary
curve with larger values of the coordinate $y$ as the upper boundary of $A$. Let $\gamma\subset A'$ be a simple essential curve and
let $\gamma^n$ be the boundary of the connected component of $A\setminus(\gamma\cup F_n(\gamma))$
adjacent to the upper boundary of $A$. This is also a simple essential curve.
Denote the operator that replaces the curve $\gamma$ by $\gamma^n$ as $\mathcal{F}_n$.
By construction, $\mathcal{F}_n(\gamma)$ has no points below $\gamma$ and, since the map $F_n$ is exact,
$\gamma\cap\mathcal{F}_n(\gamma)\ne\emptyset$. If $\mathcal{F}_n(\gamma)\cap\gamma^+\ne\emptyset$ for some $n$,
the existence of the sought connecting orbit is trivial (indeed, take $v_1\in \mathcal{F}_n(\gamma)\cap\gamma^+$ and  $v_0=F^{-1}_n(v_1)$).

Now we continue by induction. Let $m=0$, $\gamma_0=\gamma^-$ and  inductively construct a sequence of simple essential
curves $\gamma_m\subset \bar A$, such that each point on $\gamma_m$ can be reached by a trajectory which
starts on $\gamma^-$ and has the length not longer than $m$.
Suppose we have constructed such $\gamma_m$ for some $m\ge0$.
If $\mathcal{F}_n(\gamma_m)\cap\gamma^+\ne\emptyset $ for some $j$,
the inductive process is terminated as the intersection point belongs to
the sought trajectory which starts on $\gamma^-$ and finishes on $\gamma^+$.
Otherwise define $\gamma_{m+1}$ as the boundary of that connected
component of $A\setminus(\cup_n \mathcal{F}_n(\gamma_m))$ which is adjacent to the upper boundary of $A$.
Obviously, $\gamma^+$ is a simple essential curve. Since each of the curves
$\mathcal{F}_n(\gamma_m)$ intersects $\gamma_m\subset \bar A$ (by the exactness of $F_n$) and none of them, by assumption,
intersects $\gamma^+$, it follows that $\gamma_{m+1}\subset \bar A$. By construction, $\gamma_{m+1}$ has no points below $\gamma_m$.

We claim that this process will terminate after a finite number of steps because otherwise the maps $F_n$
would have a common invariant essential curve in $\bar A$. Indeed,
suppose that the process does not terminate. Recall that the inductive process produces the sequence of curves
which is ``bounded and monotone" in the sense that $\gamma_{m+1}$ lies in the region between $\gamma_m$ and $\gamma^+$ for all $m$.
Let  $U^*$ denote that connected component
of $A\setminus(\cup_{m\ge0} \gamma_m)$ which is adjacent to the upper boundary of $A$
and let $\gamma^*$ be the boundary of $U^*$.

Let us show that $F_0(U^*)=U^*$. Indeed, suppose this is not the case. Then we cannot have $U^*\subset F_0(U^*)$, as
this would imply that the set $F_0(U^*)\setminus U^*$ has non-empty interior, hence the area
of $F_0(U^*)$ will be strictly larger than the area of $U^*$ which is impossible as $F_0$ is area-preserving.
Thus, $U^*\not\subset F_0(U^*)$, which means $F_0^{-1}(U^*)\not\subset U^*$. At the same time $F_0^{-1}(U^*)\cap U^*\neq \emptyset$,
so, since $U^*$ is connected, $F_0^{-1}(U^*)$ intersects the boundary $\gamma^*$ of $U^*$, i.e.
there exists a point $p^*\in \gamma^*$ such that $F_0(p^*)\in U^*$. By the construction of $\gamma^*$,
there is a sequence of points $p_m\in\gamma_m$ which converges to $p^*$. Consequently, since $U^*$ is open,
$F_0(p_m)\in U^*$ for all sufficiently large $m$. This is impossible as
$F_0(p_m)\in F_0(\gamma_m)$, and the curve $F_0(\gamma_m)$ cannot have points above $\gamma_{m+1}$, i.e. it cannot intersect $U^*$.

Now note that since $U^*$ is invariant with respect to the twist map $F_0$ and contains an essential invariant curve of $F_0$ (the upper boundary of $A$),
the Birkhoff theory implies that the boundary $\gamma^*$ of $U^*$ is a simple essential curve, a graph of a Lipschitz function $y(\varphi)$.
Let us show that this curve is a common invariant curve for all of the maps $F_n$. Indeed, let $\gamma^*$ be not invariant with respect to
$F_n$ for some $n$. Then, since $F_n$ is exact symplectic, $F_n(\gamma^*)$ must have points on both sides from $\gamma^*$, i.e.
there exists $p^*\in\gamma^*$ such that $F_n(p^*)\in U^*$. As $\gamma^*$ lies in the closure of the curves $\gamma_m$,
there is a sequence of points $p_m\in\gamma_m$ which converges to $p^*$. Consequently,
$F_n(p_m)\in U^*$ for all sufficiently large $m$, which is impossible as
$F_n(p_m)\in F_n(\gamma_m)$, and the curve $F_n(\gamma_m)$ cannot have points in $U^*$.

We have proved that the non-existence of the trajectory of statement 1 is equivalent to the
existence of a common invariant curve of statement 2.
\hfill$\Box$

\medskip

Theorem \ref{Thm:Birkhoff} is valid for any two non-intersecting essential curves in $A'$: either they are connected by
an orbit of the iterated function system, or there is an essential curve $\gamma^*$ between them which is invariant with respect to all
maps $F_n$. It follows that the absence of a common invariant essential curve in $\bar A$ is equivalent to the existence of an orbit
of the iterated function systems which connects $int(A^+)$ with $int(A^-)$ (move the curves $\gamma^+$ and $\gamma^-$ inside $int(A^+)$ and, respectively,
$int(A^-)$, and apply Theorem \ref{Thm:Birkhoff} to these curves). Since the existence of such orbit is an open property, Theorem
\ref{Thm:Birkhoff} implies that the cylinder $\bar A$ contains no essential curve invariant with respect to all maps $F_0,\dots, F_N$ for
an {\em open} set of maps from $\mathcal V$. In the next Section we show that this set of maps is also {\em dense} in $\mathcal V$.
This will finish the proof of the Main Theorem: it follows immediately from Theorem \ref{Thm:Birkhoff} and Lemma \ref{lmsw} that for any
map $\tilde\Phi$ from this open and dense set any two neighbourhoods of $\gamma^-$ and $\gamma^+$ are connected by $\tilde\Phi$.

\section{Simultaneous destruction of all obstruction curves}

We finish the proof of the Main Theorem by showing that for a dense subset of the neighbourhood $\mathcal V$ of our original map $\Phi$
the maps $F_0,F_1,\dots,F_N$ do not have a common essential invariant curve. As we mentioned, it is enough for us to prove it for any given fixed $N$.
As $F_0$ is a twist map, we can restrict the problem to Lipshitz invariant curves only. Recall that for any map $\tilde \Phi$ from $\mathcal V$
there exists a compact normally-hyperbolic invariant cylinder $A$. We introduce coordinates $(y,\varphi)$ on $A$ such that the restriction $F_0$
of $\tilde\Phi$ on $A$ has a twist property. Namely, by denoting $F_0: (y,\varphi)\mapsto (\bar y,\bar \varphi)$ we have
$$\frac{\partial \bar y}{\partial\varphi}\neq 0$$
for all $(y,\varphi)\in A$. By a Birkhoff theorem, every essential invariant curve of $F_0$ is Lipschitz:
$$y=y(\varphi), \qquad |y(\varphi_1)-y(\varphi_2)|\leq L|\varphi_1-\varphi_2|,$$
where the Lipschitz constant $L$ satisfies
$$
L\leq \sup_{v\in\bar A} \;\max\left\{\left|\frac{\partial\bar \varphi}{\partial \varphi}\right|/\left|\frac{\partial\bar \varphi}{\partial y}\right|,\;
\left|\frac{\partial\bar y}{\partial y}\right|/\left|\frac{\partial\bar \varphi}{\partial y}\right|\right\}\,.
$$
Since all the maps from $\mathcal V$ are $C^1$-close and the corresponding cylinders $A$ are $C^1$-close as well, we can take the constant $L$ the
same for all maps from $\mathcal V$.

By assumptions of the Main Theorem and by Proposition \ref{Prop:manyb}, we have a compact subcylinder $\bar A$ in $A$ such that $N\geq 8$
scattering maps are defined on a neighbourhood $A'$ of $\bar A$. The cylinder $\bar A$ depends continuously on the map $\tilde\Phi$,
so we can choose $A'$ to be the same (in appropriately chosen coordinates $(y,\varphi)$) for all maps close to $\tilde\Phi$. We can also assume that the maps
$F_1,\dots,F_N$ are defined in some neighbourhood of the closure of $A'$. Note that the scattering maps depend continuously on the map $\tilde\Phi$ in the following
sense: if $\Phi$ is $C^2$-close to $\tilde\Phi$, then the corresponding scattering maps are $C^1$-close.

\begin{theorem}\label{Thm:nocurve} Arbitrarily close to any map $\tilde\Phi$ in $\mathcal V$ there exists a map for which
the corresponding scattering maps $F_1,\dots,F_8$ have no common $L$-Lipschitz invariant curves in $A'$.
\end{theorem}

\noindent
{\em Proof.} Consider the space of all $L$-Lipshitz functions $y=y(\varphi)$ endowed by the $C^0$-metrics. Let $\mathcal L$ be the subset
of this space which consists of all functions whose graphs lie in the closure of $A'$ and are invariant, simultaneously,
for all the scattering maps $F_1, \dots, F_8$ generated by the map $\tilde\Phi$. If $\mathcal L=\emptyset$, there is nothing to prove.
If $\mathcal L\neq \emptyset$, we note that $\mathcal L$ is compact, so given any $\delta>0$ there is a finite set of $L$-Lipshitz curves
$C_1,\dots,C_q$ such that each of them is invariant with respect to all the maps $F_1, \dots, F_8$ and every other common invariant $L$-Lipshitz curve
lies in the $\delta$-neighbourhood of one of the curves $C_s$, i.e. it belongs to the cylinder $A_s:=\{|y-y_s(\varphi)|\leq \delta\}$ where
$y=y_s(\varphi)$ is the equation of the curve $C_s$. Moreover, the set of the $L$-Lipshitz common invariant curves of the scattering maps
depends upper-semicontinuously on the map $\tilde \Phi$ (if we have a sequence of maps $\Phi^{(k)}$ that converges to $\tilde\Phi$ in $C^2$, then
the corresponding scattering maps $F^{(k)}_j$ converge to the scattering maps $F_j$ in $C^1$; and if the maps $F^{(k)}_j$ each have an $L$-Lipshitz invariant
curve, then the set of the limit points of these curves as $k\to+\infty$ is the union of a set of $L$-Lipshitz curves each of which is invariant
with respect to the scattering maps $F_j$). Thus, for all maps from $\mathcal V$ which are sufficiently close to $\tilde\Phi$, every
common invariant $L$-Lipshitz curve of the scattering maps that lies in $A'$ lies entirely in one of the cylinders $A_1,\dots,A_q$.

Below (see (\ref{ddfr1})) we will fix, once and for all, a certain value of $\delta>0$ which will give us a finite set of these cylinders $A_s$.
We will show for each such cylinder $A_s$
that arbitrarily close to $\tilde\Phi$ in $\mathcal V$ there exists a map for which
the corresponding scattering maps $F_1,\dots,F_8$ have no common $L$-Lipschitz invariant curves in $A_s$. This will prove the theorem.
Indeed, the absence of the common invariant $L$-Lipshitz curves in any given (open) cylinder is an open property. So, we first perturb the map $\tilde\Phi$
such that to kill all common invariant $L$-Lipshitz curves in the cylinder $A_1$, then we add one more small perturbation to kill all
common invariant $L$-Lipshitz curves in $A_2$ - by choosing the perturbation small enough we guarantee that no new common invariant $L$-Lipshitz curves
emerge in $A_1$, etc.. Then, after finitely many steps of the procedure we will have all the cylinders $A_1,\dots,A_q$ cleaned of the possible
common invariant $L$-Lipshitz curves.

\smallskip

Let $R>1$ be a constant that bounds the derivatives of the scattering maps:
\begin{equation}\label{rdf0}
\left\|\frac{\partial F_j}{\partial (y,\varphi)}\right\|<R
\end{equation}
for all $(y,\varphi)\in A'$, $j=1,\dots,8$, and all $\Phi\in\mathcal V$ that close enough to $\tilde\Phi$.
Recall that $\varphi$ is an angular variable that runs a circle $\mathbb S^1$; we assume that the length of the circle is $2\pi$.
Choose 4 arcs $J_i\subset {\mathbb S}^1$, $i\in\{1,2,3,4\}$,
$J_1\cup J_2=J_3\cup J_4=\mathbb S^1$. Moreover, denote $J_{ik}=J_i\setminus J_k$ and let us assume
that $J_{12}$, $J_{34}$, $J_{21}$ and $J_{43}$
are disjoint and lie in $\mathbb S^1$ in the same order as they are listed here (following the orientation on the circle
$\mathbb S^1$). Neither of the arcs $J_i$ constitutes the whole circle, so their lengths are smaller than $2\pi$. Choose any
$L$-Lipshitz curve $C: y=y_{_C}(\varphi)$ which is invariant with respect to all maps $F_1,\dots,F_8$.
Each arc $J_i$ corresponds to an arc $\hat J_i:\{y=y_{_C}(\varphi),\varphi\in J_i\}$ of the curve $C$.
Since $C$ is invariant with respect to each of the maps $F_j$, the image $F_j(\hat J_i)$ also lies in $C$. Hence it
is given by $F_j(\hat J_i): \{y=y_{_C}(\varphi),\varphi\in {\bar J}_i^{j}\}$ where ${\bar J}_i^{j}$ is a certain arc in $\mathbb S^1$
which does not cover the whole of $\mathbb S^1$, so its length is strictly less than $2\pi$. Since the set $\mathcal L$
of all common invariant $L$-Lipshitz curves is compact, we have
\begin{equation}\label{kdfk}
K=\max_{C\in\mathcal L} \; \max_{i,j} \; {\rm length}({\bar J}_i^j)<2\pi.
\end{equation}
Now, we choose
\begin{equation}\label{ddfr1}
\delta=\frac{2\pi-K}{R}>0.
\end{equation}

As it was explained above, the compactness of $\mathcal L$ implies that every possible common invariant $L$-Lipshitz curve lies in
one of a finitely many cylinders $A_s$;
each of these cylinders is the $\delta$-neighbourhood of some invariant $L$-Lipshitz curve $C_s:\{y=y_s(\varphi)\}$.
Take any of these cylinders. Note that, by virtue of (\ref{rdf0}), the image $F_j(A_s\cap\{\varphi\in J_i\})$ lies inside the $(R\delta)$-neighbourhood
of the curve $F_j(C_s\cap\{\varphi\in J_i\})$. This curve is a subset of the invariant curve $C_s$, and it corresponds to an interval of $\varphi$ values
such that the length of this interval does not exceed the constant $K$ defined by (\ref{kdfk}). Thus, by (\ref{ddfr1}),
\begin{equation}\label{idns}
F_j(A_s\cap\{\varphi\in J_i\})\subset \{|y-y_s(\varphi)|<R\delta,\; \varphi\in \hat J_{sij}\}
\end{equation}
where $\hat J_{sij}$ is a certain arc whose length is strictly less than $2\pi$, i.e. it does not cover the entire $\mathbb S^1$.
As $F_j$ depends continuously on the map $\tilde\Phi$, inclusion (\ref{idns}) holds for all maps from $\mathcal V$ which are close enough to
$\tilde\Phi$.

Now, let us imbed the map $\tilde \Phi$ into a two-parameter analytic family of maps $\Phi_{\mu_1,\mu_2}$
from $\mathcal V$ such that $\Phi_0=\tilde\Phi$.
We will show (Lemmas \ref{rlipmov} and \ref{hampr}) that this family can be chosen such that there exist arbitrarily small values of $\mu=(\mu_1,\mu_2)$ for which
the scattering maps $F_1,\dots,F_8$ defined by the map $\Phi_\mu$ have no common $L$-Lipschitz invariant curves in the cylinder $A_s$.
The map $\Phi_\mu$ that corresponds to a small value of $\mu$ is a small perturbation of $\tilde\Phi$, so this gives us the required
arbitrarily small perturbations that clear the cylinder $A_s$ of the common $L$-Lipshitz invariant curves of the scattering maps.
By performing this perturbations consecutively for each of the cylinders $A_1,\dots, A_q$ we will obtain the result of the theorem.

Note that the invariant cylinder $A$, its stable and unstable manifolds, as well as the strong stable and strong unstable foliations
depend smoothly on $\mu$, therefore the scattering maps also depend smoothly on $\mu$.
Denote
$$F_j: (y,\varphi)\mapsto (Y_j(y,\varphi,\mu), \Psi_j(y,\varphi,\mu)).$$
Let our family $\Phi_\mu$ be chosen such that for all $(\varphi,y)\in A_s$
\begin{equation}\label{kpm}\!\!\!\!\!\!\!\!\!\!\!\!\!
\left\|\frac{\partial \Psi_j}{\partial (\mu_1,\mu_2)}\right\|<1 \qquad\mbox{for all $j=1,\dots,8$},
\end{equation}
\begin{equation}\label{yyyy}
\!\!\!\!\!\!\!\!\!\!\!\!\!\!\!\!\!\!\!\!\!\!\!\!\!\!\!\!
\left|\frac{\partial Y_{1,2,3,4}}{\partial \mu_2}\right| <1, \qquad \left|\frac{\partial Y_{5,6,7,8}}{\partial \mu_1}\right|< 1,
\end{equation}
\begin{equation}\label{jjj}
\!\!\!j=1,2: \qquad
\frac{\partial Y_{j}}{\partial \mu_1}>\;\;2(L+1) \;\;\mbox{ and }\;\; \frac{\partial Y_{j+4}}{\partial \mu_2} > \;\;2(L+1)
\qquad\mbox{when  } \Phi_j(\varphi,y,\mu)\in J_j,
\end{equation}
\begin{equation}\label{jjj-}
\!\!\!
j=3,4:\qquad \frac{\partial Y_{j}}{\partial \mu_1}< -2(L+1) \;\;\mbox{ and }\;\; \frac{\partial Y_{j+4}}{\partial \mu_2} < -2(L+1)
\;\quad\mbox{when } \Phi_j(\varphi,y,\mu)\in J_j,
\end{equation}
where $L$ is the Lipschitz constant in the condition of the theorem, and $J_{1234}$ are the arcs defined above.

We will show below (Lemma \ref{hampr}) that such family $\Phi_\mu$ exists. Thus, we will prove the theorem by showing the following result.

\begin{lemma}\label{rlipmov}
For every family of maps $\Phi_\mu$, $\mu=(\mu_1,\mu_2)$, such that the derivatives of the scattering maps $F_1,dots,F_8$ satisfy estimates (\ref{kpm})-(\ref{jjj-})
for all $(\varphi,y)\in A_s$, the set of parameter values for which the scattering maps $F_1,dots,F_8$ have an $L$-Lipshitz common invariant essential
curve in $A_s$ has measure zero. In particular, there exist arbitrarily small values of $\mu$ for which the maps $F_1,dots,F_8$ have no
$L$-Lipshitz common invariant essential curves in the cylinder $A_s$.
\end{lemma}

{\em Proof.} Take any two, may be equal, values of $\mu$: $\mu=\mu^*$ and $\mu=\mu^{**}$, such that at $\mu=\mu^*$
the maps $F_1,\dots,F_8$ have a common $L$-Lipschitz invariant curve ${\mathcal L}^*: \{y=y^*(\varphi), \varphi\in \mathbb S^1\}\subset A_s$ and
at $\mu=\mu^{**}$ they have a common $L$-Lipschitz invariant curve ${\mathcal L}^{**}: \{y=y^{**}(\varphi), \varphi\in \mathbb S^1\}\subset A_s$.
Let us show that the following condition holds:
\begin{equation}\label{qqq}
\|\mu^*-\mu^{**}\|\leq R |y^*(0)-y^{**}(0)|,
\end{equation}
where $R$ is defined in (\ref{rdf0}) and $\|\mu\|=\max\{|\mu_1|,|\mu_2|\}$.

We note that without loosing in generality we may assume that
\begin{equation}\label{yps}
y^*(0)\geq y^{**}(0),
\end{equation}
\begin{equation}\label{m21ll}
|\mu_2^*-\mu_2^{**}|\leq |\mu_1^*-\mu_1^{**}|\qquad\mbox{and}\qquad  \mu_1^*\geq \mu_1^{**}.
\end{equation}
If necessary, these inequalities can be achieved by swapping $y$ and $(-y)$,  $\mu$ and $(-\mu)$,  $F_1\leftrightarrow F_3$,
$F_2\leftrightarrow F_4$, $F_5\leftrightarrow F_7$, $F_6\leftrightarrow F_8$, as well
as $\mu_1\leftrightarrow \mu_2$ and $F_{1,2,3,4}\leftrightarrow F_{5,6,7,8}$.
Conditions (\ref{kpm})-(\ref{jjj-}) are symmetric with respect to these changes.

Now suppose (\ref{qqq}) is not true, i.e.
\begin{equation}\label{nolmy}
0\leq y^*(0)-y^{**}(0)<\frac{\Delta\mu}{R},
\end{equation}
where
$$\Delta\mu=\mu_1^*-\mu_1^{**}>0.$$
Since $J_3\cup J_4=\mathbb S^1$, we have that $\varphi=0$ lies in least in one of the arcs $J_3$ or $J_4$.
For definiteness, we assume $0\in J_3$. Let $(\varphi^*,\bar y^*)=F_3(0,y^*(0),\mu^*)$  and $(\varphi^{**},\bar y^{**})=F_3(0,y^{**}(0),\mu^{**})$,
i.e.
$$
\begin{array}{rcllcl}
\varphi^*&=&\Psi_3(0,y^*(0),\mu^*), \qquad&\bar y^*&=&Y_3(0,y^*(0),\mu^*),\\
\varphi^{**}&=&\Psi_3(0,y^{**}(0),\mu^{**}),\qquad &\bar y^{**}&=&Y_3(0,y^{**}(0),\mu^{**}).
\end{array}
$$
Standard estimates based on the mean value theorem and formulas
(\ref{rdf0}),(\ref{kpm}),(\ref{yyyy}),(\ref{jjj-}),(\ref{m21ll}),(\ref{nolmy})
imply that
$$|\varphi^{**}-\varphi^*| <  2\Delta\mu, \qquad \bar y^*-\bar y^{**} < -2L \Delta\mu.
$$
Since the curves $y=y^*(\varphi)$ and $y=y^{**}(\varphi)$ are invariant with respect to $F_3$
(at $\mu=\mu^*$ and $\mu=\mu^{**}$ respectively),
it follows that $\bar y^*=y^*(\varphi^*)$, $\bar y^{**}=y^{**}(\varphi^{**})$. Because of the $L$-Lipschitz property,
we find that
$$y^*(\varphi^*)-y^{**}(\varphi^*)=\bar y^*-\bar y^{**}+y^{**}(\varphi^{**})-y^{**}(\varphi^*)<-2L\Delta\mu+2L\Delta\mu <0.$$
Then taking into account (\ref{yps}) we conclude that
$${\mathcal L}^*\cap{\mathcal L}^{**}\ne \emptyset.$$

Now, let us call an arc $I\subset \mathbb S^1$ of the $\varphi$ values {\em positive} if
$y^*(\varphi)>y^{**}(\varphi)$ for all $\varphi\in int(I)$,
and $y^*(\varphi)=y^{**}(\varphi)$ at the end points of $I$.
We call an arc {\em negative}, if $y^*(\varphi)=y^{**}(\varphi)$ at its end points and
$y^*(\varphi)<y^{**}(\varphi)$ on its interior.
It is convenient to allow arcs to have empty interiors, i.e. any point from
${\mathcal L}^*\cap{\mathcal L}^{**}$ is considered to be both a positive and a negative arc
at the same time.

\smallskip

We have just proved that there is at least one negative and at least one positive arc.
Let $I$ be a positive arc. Let ${\mathcal L}^*_I=\{\,y=y^*(\varphi),\varphi\in I\,\}$ and ${\mathcal L}^{**}_I=\{\,y=y^{**}(\varphi),\varphi\in I\,\}$ be
the corresponding pieces of the curves ${\mathcal L}^*$ and ${\mathcal L}^{**}$, and let ${\mathcal D}_I=\{\,y^*(\varphi) \geq y\geq y^{**}(\varphi),\varphi\in I\,\}$
be the region bounded by ${\mathcal L}^*_I$ and ${\mathcal L}^{**}_I$.
Let us show that if $I\subseteq J_j$ for $j=1$ or $j=2$, then, with this $j$, the image of ${\mathcal L}^*_I$ by the map $F_j$ at
$\mu=\mu^*$ lies strictly
inside ${\mathcal L}^*_{I'}$ which corresponds to a positive arc $I'$ and
\begin{equation}\label{lenpp}
\mathrm{length}(I')>\Delta\mu>0,
\end{equation}
\begin{equation}\label{arspp}
\mathrm{area}({\mathcal D}_{I'})>\mathrm{area}({\mathcal D}_{I}).
\end{equation}

Indeed, denote as $F_j^*$ the map $F_j$ at $\mu=\mu^*$ and $F_j^{**}$ the map $F_j$ at $\mu=\mu^{**}$. Take any point
$M(\varphi,y^*(\varphi))\in {\mathcal L}^*_I$, so $\varphi\in I$. Let $M^*(\varphi^*,y^*(\varphi^*))\in {\mathcal L}^*$ be the image
of $M$ by the map $F_j^*$, and $M'(\varphi',y')\in F_j^{**}({\mathcal L}^*_I)$ be the image
of $M$ by the map $F_j^{**}$. Since $I$ is a positive arc, we have that for any $\varphi\in I$
the point $M$ is either on
the curve ${\mathcal L}^{**}$ or above it. Since ${\mathcal L}^{**}$ is invariant with respect to $F_j^{**}$,
the point $M'=(\varphi',y')$ also does not lie below ${\mathcal L}^{**}$, i.e.
\begin{equation}\label{yyp}
y'\geq y^{**}(\varphi').
\end{equation}
We have
$$\varphi'=\Psi_j(\varphi,y^*(\varphi),\mu^{**}),\qquad  y'=Y_j(\phi,y^*(\phi),\mu^{**}),$$
$$\varphi^*=\Psi_j(\varphi,y^*(\varphi),\mu^*),\qquad y^*(\varphi^*)=Y_j(\varphi,y^*(\varphi),\mu^*).$$
Thus, inequalities (\ref{kpm})-(\ref{jjj}) imply that
$$|\varphi^*-\varphi'|<\Delta\mu,\qquad y^*(\varphi^*)-y'>(2L+1)\Delta\mu$$
(recall that we assume $I\subseteq J_j$, hence $\varphi\in J_j$). By (\ref{yyp}) and the $L$-Lipschitz property
of ${\mathcal L}^{**}$ we obtain
\begin{equation}\label{posiyp}
y^*(\varphi^*)-y^{**}(\varphi^*)>(L+1)\Delta\mu>0,
\end{equation}
and
\begin{equation}\label{yposip}
y^*(\varphi')-y'>(L+1)\Delta\mu>0.
\end{equation}

Denote $\tilde F_j(\varphi)=\Psi_j(\varphi,y^*(\varphi),\mu^*)$, i.e., $\tilde F_j$ is the restriction
of the map $F_j^*$ on the invariant curve ${\mathcal L}^*$.
We have just showed that if $\varphi\in I$, where $I\subseteq J_j$ is a positive arc,
then $\varphi^*=\tilde F_j(\varphi)$ satisfies (\ref{posiyp}), i.e. it is inside some positive arc $I'$.
Moreover, at the end points of $I'$ we must have $y^*-y^{**}=0$ while
at the points of $\tilde F_j(I)\Subset I'$ we have $y^*-y^{**}>L\Delta\mu$ by (\ref{posiyp}), hence
the length of $I'$ is bounded from below as in (\ref{lenpp}), by virtue of the
$2L$-Lipschitz property of the function $y^*(\varphi)-y^{**}(\varphi)$.

We have shown that $F_j^*(\mathcal L^*_I)\subset \mathcal L^*_{I'}$ and $F_j^{**}(\mathcal L^{**}_I)\subset \mathcal L^{**}_{I'}$
where $I'$ is a positive arc. As the point $M$ runs $\mathcal L^*_I$, the point $M'(\varphi',y')$ runs the curve $\mathcal L'=F^{**}(\mathcal L^*_I)$,
and it follows from (\ref{yyp}),(\ref{yposip}) that the curve $\mathcal L'$ lies between $\mathcal L^*$ and $\mathcal L^{**}$, strictly below
$\mathcal L^*$. Since the end points of $\mathcal L'$ coincide with the end points of $F_j^{**}(\mathcal L^{**}_I)$
and the latter lie inside $\mathcal L^{**}_{I'}$, it follows that $\mathcal L'$ lies between $\mathcal L^*_{I'}$ and $\mathcal L^{**}_{I'}$, strictly
below $\mathcal L^*_{I'}$. Therefore the area of the region $F^{**}_j(\mathcal D_I)$ bounded by the curves $\mathcal L'$ and $F_j^{**}(\mathcal L^{**}_I)$
is strictly smaller than the area of the region $\mathcal D_{I'}$ bounded by the curves $\mathcal L^*_{I'}$ and $\mathcal L^{**}_{I'}$. As the map $F_j$
is area-preserving, $\mathrm{area}(F^{**}_j\mathcal D_I)=\mathrm{area}(\mathcal D_I)$, and (\ref{arspp}) follows.

Thus, we start with any positive arc $I$ which is contained entirely inside $J_1$ or $J_2$ and obtain a sequence
$I_s$ of positive arcs such that $I_0=I$ and $\tilde F_{j_s}(I_s)\subset I_{s+1}$, where we chose $j_s=1$ if $I_s\subseteq J_1$, and
$j_s=2$ if $I_s\subseteq J_2$ and $I_s\not\subseteq J_1$. If for some $s$ the arc $I_s$ is not entirely contained neither in $J_1$
nor in $J_2$, the sequence is terminated. By (\ref{arspp}), the area of the region $\mathcal D_{I_s}$ is a strictly increasing function of $s$,
so the arcs with different $s$ can never coincide. The definition of a positive arc implies that the intersection of interiors for two
different positive arcs is always empty. Thus, the arcs $int(I_s)$ are mutually disjoint. By (\ref{lenpp}), no more than $\frac{2\pi}{\Delta\mu}$
of such arcs can coexist in $\mathbb S^1$. We conclude
that the sequence $I_s$ must terminate. This means the last arc in the sequence is not contained entirely neither in
$J_1$ nor in $J_2$, i.e. we have proved that there is a positive arc $I^+$ such that both $I^+\cap J_{12}\ne\emptyset$ and $I^+\cap J_{21}\ne\emptyset$.

Similarly, one proves that there exists a negative arc $I^-$ such that
$I^-\cap J_{34}\ne\emptyset$ and $I^-\cap J_{43}\ne\emptyset$. Since $J_{12}$, $J_{34}$, $J_{21}$ and $J_{43}$ are placed on $\mathbb S^1$
in this order, we find that the interiors of $I^+$ and $I^-$ intersect, which is impossible by the definition of
positive and negative arcs. Thus, by contradiction, we have  established estimate (\ref{qqq}).

\smallskip

Let $\mathcal M\subset \mathbb R^2$ be the set of all $\mu$ such that the maps $F_1,\dots,F_8$ have at least one common  $L$-Lipschitz invariant curve
in the cylinder $A_s$. Let $\mathcal Y$ be the set which consists of all intersection points of these curves with the axis $\varphi=0$.
By (\ref{qqq}), for each $y_0\in{\mathcal Y}$ there is exactly one $\mu \in \mathcal M$
such that the corresponding system of scattering maps has a common $L$-Lipshitz invariant curve that lies in $A_s$ and intersects the line $\varphi=0$ at $y=y_0$.
Estimate (\ref{qqq}) also implies  that $y_0\mapsto \mu$ is an $R$-Lipschitz function $\mathcal Y \to \mathcal M$.
For any Lipschitz function from a subset of $\mathbb R$ to $\mathbb R^2$, the Lebesgue measure
of the image vanishes. Thus, as $\mathcal Y$ is a subset of an interval, it follows that  $\mathrm{mes}(\mathcal M)=0$, which gives the lemma.
\hfill$\Box$\\

We stress that Lemma \ref{rlipmov} holds for any family of symplectic maps $\Phi_\mu$ for which conditions (\ref{kpm})-(\ref{jjj-}) are satisfied.
In order to finish the prove of the theorem, it remains to show that such family can be taken inside the space $\mathcal V$ of analytic
exact-symplectic maps. This is given by the lemma below.

\begin{lemma}\label{hampr}
Any map $\tilde\Phi\in \mathcal V$ can be imbedded into
an analytic family of analytic exact-symplectic maps $\Phi_\mu$ that satisfies conditions (\ref{kpm})-(\ref{jjj-}).
\end{lemma}

{\em Proof.} We define $\Phi_\mu=X_\mu\circ \tilde\Phi$, where $X_\mu$ is an analytic family of exact-symplectic maps
such that $X_0=id$. Specifically, we set $X_mu=X^{(1)}_{\mu_1}\circ X^{(2)}_{\mu_2}$ where
$X^{(i)}_{\mu_i}$ is the time-$\mu_i$
shift by the orbits of the vector field defined by an analytic Hamiltonian function $H_i$ ($i=1,2$).
Since we are interested in small $\mu$, it is enough to check the fulfillment of conditions
(\ref{kpm})-(\ref{jjj-}) at $\mu=0$ only. Thus the family $\Phi_\mu=X_\mu\circ \tilde\Phi$ will satisfy
(\ref{kpm})-(\ref{jjj-}) for all small $\mu$, provided the conditions
\begin{equation}\label{kpmi}\begin{array}{ll}&\displaystyle
\left|\frac{\partial \Psi_j}{\partial \mu_i}\right|_{\mu_i=0}<1 \qquad (j=1,\dots,8),\\ \mbox{for all } \varphi\in \mathbb S^1 :\quad\\ &\displaystyle
\left|\frac{\partial Y_j}{\partial \mu_i}\right|_{\mu_i=0}<1 \qquad
(j=9-4i,\dots,12-4i),
\end{array}
\end{equation}
\begin{equation}\label{j3i-}\!\!\!\!\!\!
\begin{array}{l}\displaystyle
\mbox{for all $\varphi\in J_j$ with $j=1,2$:} \qquad
\left.\frac{\partial Y_{j+4(i-1)}}{\partial \mu_i}\right|_{\mu_i=0}>\;\;2(L+1),\\ \\ \displaystyle
\mbox{for all $\varphi\in J_j$ with  $j=3,4$:}\qquad
\left.\frac{\partial Y_{j+4(i-1)}}{\partial \mu_i}\right|_{\mu_i=0}< -2(L+1),
\end{array}
\end{equation}
are satisfied by the scattering maps for the families $\Phi_{\mu_i}^{(i)}=X^{(i)}_{\mu_i}\circ \tilde\Phi$, $i=1,2$, for all $(\varphi,y)\in A_s$.

Let us build a family of maps $X^{(1)}_{\mu_{_1}}$ for which these conditions are satisfied (the construction for $i=2$ is the same).
Conditions (\ref{kpmi}),(\ref{j3i-})
are strict and involve only the first derivatives of the scattering maps.
A $C^2$-small change of the family $\Phi_{\mu_{_1}}^{(1)}$ leads to a $C^1$-small change of the strong-stable and strong-unstable
foliations and, therefore, a $C^1$-small change of the scattering maps. Thus, it is enough
to build a $C^2$-smooth family of maps $X^{(1)}_{\mu_{_1}}$ (generated by a $C^3$-smooth Hamiltonian $H^{(1)}$)
such that the corresponding scattering maps satisfy (\ref{kpmi})-(\ref{j3i-}). Then for any sufficiently
$C^3$-close approximation of $H^{(1)}$ by an analytic Hamiltonian (the analiticity of $H^{(1)}$ and $H^{(2)}$
is needed for the family $\Phi_\mu$ to be analytic, i.e. lie in $\mathcal V$) conditions (\ref{kpm})-(\ref{jjj-}) will still be satisfied.

We will make the $C^3$-smooth Hamiltonian $H^{(1)}$ localised in a small neighbourhood of the cylinders
$\tilde\Phi(B_1)$, $\tilde\Phi(B_2)$, $\tilde\Phi(B_3)$, $\tilde\Phi(B_4)$. Thus, the maps $X^{(1)}_{\mu_{_1}}$
will differ from identity only in a small neighbourhood of these cylinders, so the maps
$\Phi^{(1)}_{\mu_{_1}}$ will differ from $\tilde\Phi$ in a small neighbourhood of the cylinders $B_1,\dots,B_4$
only. The perturbation we build near one of these cylinders does not affect the scattering maps near the other
cylinders, so we restrict our attention to the cylinder $B_1$ only. We further omit the index "1" whenever possible and
denote $\tau=\mu_1$. Thus we consider a homoclinic cylinder $B$ and continue with building a $C^3$-smooth
Hamiltonian $H$ localised in a small neighbourhood of the cylinder $\tilde\Phi(B)$ such that for
the corresponding flow map $X_\tau$ the derivative with respect to $\tau$ of the scattering map $F$ defined by the map $\Phi_\tau=X_\tau\circ \tilde\Phi$
satisfies, for all $(\varphi,y)\in A_s$, the following inequalities:
\begin{equation}\label{kpmi1}\begin{array}{l} \displaystyle
\left|\frac{\partial \Psi}{\partial \tau}\right|_{\tau=0}<1 \;\; \mbox{ for all } \varphi\in\mathbb S^1,\\ \\ \displaystyle
\left.\frac{\partial Y}{\partial \tau}\right|_{\tau=0} >2(L+1)\;\; \mbox{ for all } \varphi\in J,
\end{array}
\end{equation}
where $J$ is a certain arc that does not contain the whole $\mathbb S^1$, and
\begin{equation}\label{idm1s}
F(A_s\cap \{\varphi\in J\}) \subset \{\varphi\in \hat J\}
\end{equation}
where $\hat J$ is an arc that does not contain the whole of $\mathbb S^1$ (see (\ref{idns})).

Let $w^u$ denote a piece of the unstable manifold $W^u(A)$ that contains the cylinder $B$ (i.e. $w^u$ is a small neighbourhood of the
cylinder $B$ in $W^u(A)$) and $w^s$ be a small neighbourhood of $\tilde\Phi(B)$ in $W^s(A)$, so $B=\tilde\Phi(w^u)\cap w^s$.
Since the map $\Phi_\tau$ differs from $\Phi$ in a small neighbourhood of the cylinder $B$ only, the pieces $w^u$ and $w^s$ do not
depend on $\tau$, nor the strong unstable foliation of the piece of $W^u(A)$ between $A$ and $w^u$ depends on $\tau$, neither
the strong stable foliation of the piece of $W^s(A)$ between $w^s$ and $A$ does. Thus, given any $C^1$-family of cylinders $B_\tau$
close to $B$ the projection map $\pi^u_{B_\tau}:B_{\tau}\to A$ by the leaves of the strong unstable foliation is of class $C^1$;
moreover, if two such families of cylinders are $C^1$-close, then the corresponding projection maps $\pi^u_{B_\tau}$ are also $C^1$-close.
The same holds true for the projection map $\pi^s_{B_{\tau}'}:B_{\tau}'\to A$ by the leaves of the strong stable foliation, where
we denote as $B'_\tau$ any $C^1$-family of cylinders close to $\tilde\Phi(B)$. As the perturbation $X_\tau$ is localised in a small neighbourhood of the
cylinder $\tilde\Phi(B)$, we find that the scattering map $F$ satisfies
\begin{equation}\label{pffms}
F=F_0^{-1}\circ\pi^s_{B'_\tau}\circ X_\tau\circ\tilde\Phi\circ(\pi^u_{B_\tau})^{-1},
\end{equation}
where $B_\tau=w^u\cap \Phi_{\tau}^{-1}(w^s)$ is a homocinic cylinder close to $B$, and $B'_{\tau}=\Phi_\tau(B_\tau)$. If we add to the family $X_\tau$
any $C^1$-small perturbation localised in a small neighbourhood
of $\tilde\Phi(B)$, this will result in $C^1$-small perturbations of the family of cylinders $B'_\tau$ and $B_\tau$. Thus, the perturbation
to the corresponding family of scattering maps defined by (\ref{pffms}) will be also $C^1$-small. It follows that it is enough to build
a $C^1$-family of maps $X_\tau$ (generated by a $C^2$-smooth Hamiltonian $H$) localised in a small neighbourhood of the cylinder
$\tilde\Phi(B)$ such that  the corresponding family of scattering maps satisfies (\ref{kpmi1}). Then any $C^3$-Hamiltonian which is $C^2$-close to $H$
and is localised in a small neighbourhood of $\tilde\Phi(B)$ will produce a family of scattering maps that still satisfies (\ref{kpmi1}).

This reduction of smoothness requirement (from $H\in C^3$ to $H\in C^2$) is important since it allows to construct the Hamiltonian $H$ such
that the vector field it generates is tangent to the given homoclinic cylinder $B$ (for which only $C^2$-smoothness can be guaranteed by our
spectral gap assumptions). Once this is done, the cylinder $\tilde\Phi(B)$ will be invariant with respect to the map $X_\tau$, i.e. $\Phi_\tau(B)=\tilde\Phi(B)$
for all $\tau$. This means the trajectory of $B$ remains the same for all $\tau$, i.e. it remains a homoclinic cylinder. Thus, formula (\ref{pffms})
for the scattering map will recast as
\begin{equation}\label{pffms0}
F=F_0^{-1}\circ\pi^s_{\tilde\Phi(B)}\circ X_\tau\circ\tilde\Phi\circ(\pi^u_{B})^{-1},
\end{equation}
and the only $\tau$-dependent term in the right-hand side is $X_\tau$.

In order to build the required Hamiltonian, we introduce $C^2$-coordinates $(x,v)$ near $\tilde\Phi(B)$ such that
the cylinder $\tilde\Phi(B)$ is given by $x=0$ (so $v$ gives the coordinates on the cylinder and $x$ runs a neighbourhood of zero in $\mathbb R^{2d-2}$.
The cylinder is transverse to the strong-stable and strong-unstable foliations, so
if we denote as $N(v)$ the direct sum of the tangents to the leaves of the strong-stable and unstable foliations that pass through the point
$(x=0,v)\in \tilde\Phi(B)$, then the field $N(v)$ will have a form $dv=P(v)dx$. Note that $N$ depends smoothly on $v$ (as the fields of tangents to the strong
stable and strong unstable leaves are smooth when the large spectral gap assumption (\ref{nhm}) is fulfilled), i.e. the function $P(v)$ is at least $C^1$.
As the homoclinic cylinder $\tilde\Phi(B)$ belongs both to the stable and unstable manifolds of $A$, it follows from Proposition \ref{Pro:ort} that
a vector is tangent to $\tilde\Phi(B)$ if an only if it is $\Omega$-orthogonal to $N$. Thus, the vector field $\tilde X=\Omega^{-1}\nabla H$ generated
by the Hamiltonian $H$ will be tangent to $\tilde\Phi(B)$ if the gradient of $H$ is orthogonal to $N$ at the points of $\tilde\Phi(B)$, i.e.
\begin{equation}\label{hxvq}
\frac{\partial H}{\partial x}(0,v)+\frac{\partial H}{\partial v}(0,v)P(v)=0.
\end{equation}
This condition is satisfied e.g. by any function of the form
$$H(x,v)=h(v)-\sum_{i=1}^{2d-2} x_i\int p_i(v_1+s_1x_i,v_2+s_2x_i)\xi(s_1,s_2) d^2s$$
where $h$ is any $C^2$-function on $\tilde\Phi(B)$, the vector-function $p(v)=(p_1(v),\dots,p_{2d-2}(v))$ is given by $p(v)=h'(v)P(v)$, the
$x_i$'s are the coordinates of the vector $x$, and $(v_1,v_2)=v$, and
$\xi$ is a $C^2$-smooth function on a plane, localised in a small neighbourhood of zero, such
that $\int \xi(s)d^2s=1$. Integrating by parts, we find
$$\frac{\partial H}{\partial x_i}=\int p_i(v+sx_i)[s\xi'(s)+\xi(s)] d^2s,\qquad
\frac{\partial H}{\partial v_j}=\frac{\partial h}{\partial v_j}(v) + \sum_{i=1}^{2d-2} \int p_i(v+sx_i)\frac{\partial x_i}{\partial s_j} d^2s.$$
By plugging $x=0$ in these formulas, we see that (\ref{hxvq}) is satisfied indeed. Since $q\in C^1$ and $\xi\in C^2$, it follows that $H\in C^2$,
so given any $C^2$-function $h$ on the cylinder $\tilde\Phi(B)$ we can extend it to a $C^2$-function $H$ defined in a neighbourhood of this cylinder,
such that the vector field generated by the Hamiltonian $H$ will be tangent to the cylinder.

As we explained above, under this condition the scattering map is given by (\ref{pffms0}), so the vector field
$$\tilde F=\left(\tilde \Psi=\left.\frac{\partial \Psi}{\partial \tau}\right|_{\tau=0}, \tilde Y=\left.\frac{\partial Y}{\partial \tau}\right|_{\tau=0}\right)$$
of the $\tau$-derivatives of the scattering map $F$ on the cylinder $A$ is given by
\begin{equation}\label{pft}
\tilde F=\frac{\partial}{\partial v} \left(F_0^{-1}\circ\pi^s_{\tilde\Phi(B)}\right)\circ \tilde X\circ\tilde\Phi\circ(\pi^u_{B})^{-1},
\end{equation}
where $\tilde X=\Omega^{-1}(v) h'(v)$ is the vector field of the flow on the cylinder $\tilde\Phi(B)$, which is generated by the Hamiltonian $h$; we denote
as $\Omega(v)$ the antisymmetric $(2\times 2)$-matrix that defines the restriction of the symplectic form on the cylinder at the point $v$. In order to
satisfy (\ref{kpmi1}), we need to have
\begin{equation}\label{kpmi2}\begin{array}{l}
|\tilde\Psi|<1 \;\; \mbox{ for all } \varphi\in\mathbb S^1,\\ \\
\tilde Y >2(L+1)\;\; \mbox{ for all } \varphi\in J.
\end{array}
\end{equation}
It is seen from (\ref{pft}) that if conditions (\ref{kpmi2}) are satisfied by $\tilde F$
for some choice of the vector field $\tilde X$, they will be satisfied by $\tilde F$
for any $C^0$-small perturbation of $\tilde X$. Thus, it is enough to find any $C^1$-smooth Hamiltonian function $h(v)$ such that
the field $\tilde F$ defined by (\ref{pft}) satisfies (\ref{kpmi2}), then for any $C^2$-smooth function which is $C^1$-close to $h$
the derivative of the scattering map $F$ with respect to $\tau$ will satisfy (\ref{kpmi1}), and the lemma will be proven.

In order to build the sought $C^1$-function $h(v)$, we introduce $C^1$-coordinates $v=(\varphi,y)$ on the cylinder $\tilde\Phi(B)$ such that the
diffeomorphism $F_0^{-1}\circ\pi^s_{\tilde\Phi(B)}: \tilde\Phi(B)\to A$ is identity. Then (\ref{pft}) recasts as
$$\tilde F=\tilde X \circ \left.F\right|_{\tau=0}$$
(see (\ref{pffms0})). As $\tilde X$ is a Hamiltonian vector field, its $\varphi$-component is given by $\displaystyle-\omega^{-1}\frac{\partial h}{\partial y}$
and the $y$-component is $\displaystyle\omega^{-1}\frac{\partial h}{\partial \varphi}$, where the $C^0$-function $\omega(\varphi,y)> 0$ is such that
$\omega(\varphi,y)\; dy\bigwedge d\varphi$ is the symplectic  form on the cylinder $\tilde\Phi(B)$. Thus, conditions (\ref{kpmi2}) take the form
$$\left|\frac{\partial h}{\partial y}\right|<\omega(\varphi,y) \;\; \mbox{ for all } (\varphi,y)\in F(A_s),$$
$$\frac{\partial h}{\partial \varphi} >2(L+1)\omega(\varphi,y)\;\; \mbox{ for all } (\varphi,y)\in F(A_s \cap\{\varphi\in J\}).$$
We finish the proof of the lemma and the theorem by noticing that these conditions are satisfied by a $y$-independent function $h$ such that
$$h(\varphi)=M \varphi \mbox{  at  } \varphi\in\hat J$$
where the constant $M$ is given by $\displaystyle M=1+2(L+1)\sup_{F(A_s)} \omega$, and the arc $\hat J$ is defined by (\ref{idm1s});
since $h$ must be periodic in $\varphi$, it is important that $\hat J$ does not cover the whole of $\mathbb S^1$. \hfill$\Box$

\end{document}